\documentclass[12pt]{article}
\usepackage{a4}

\usepackage{epsf}

\newtheorem{picbox}{Box}

\begin{document}

\centerline{\large \bf What is Aperiodic Order?}

\bigskip
\bigskip

\centerline{Michael Baake, Uwe Grimm, Robert V.\ Moody}

\bigskip
\bigskip
\bigskip

\section{Introduction}

Surely one of the most miraculous aspects of Nature is its self-organizing
ability of creating solid substances with corresponding well-defined
macroscopic properties (namely material objects of the world around us) using
vast numbers of sub-microscopic building blocks (namely atoms and molecules).
Underlying this is the mystery of long-range order. Even putting aside the
difficult kinematic questions about crystal growth, there remains a host of
profound geometric problems: what do we mean by long-range order, how is it
characterized, and how can we model it mathematically?

In crystals, like ice, sugar, and salt, many of the extraordinarily exact
macroscopic features derive from a very simple geometric idea: the endless
repetition of a (relatively) small pattern.  A small arrangement of atoms
forms a fundamental cell that constitutes a building block, copies of which
are stacked together like bricks to fill out space by periodic repetition.
Simple as this model is, it is still difficult to analyze in full mathematical
detail: there are 230 possible symmetry classes (called space groups)
theoretically available for such periodic cell arrangements, each of which is
now also known to actually exist in Nature.  However, it took almost 100 years
from the theoretical classification of the 230 space groups to the
experimental discovery of the last examples.  Nonetheless, the underlying
feature of all crystals, which appear ubiquitously in the natural world, is
their pure periodic structure in three independent directions --- their
so-called lattice symmetry.  The interesting thing is that there is striking
long-range order in Nature that does not fit into this scheme, and one
important example of this has only been discovered recently.

Early in the last century, the wonderful tool of $X$-ray diffraction was
introduced, based on much older ideas of optical scattering (which is what we
will use to explain its essence).  Initially, diffraction pictures provided
powerful evidence of the truth of the atomic theory of matter.  Over the
years, they have become a standard tool for analyzing crystals, and to detect
long-rang order through the appearance of sharp reflection spots in the
diffraction image.  The basic idea can be visualized with an optical bench
which is driven by a small laser as source for the coherent light
(Box~\ref{laser}), see \cite{atlas} for details on this, with many instructive
examples.

Diffraction pictures of crystals display beautiful point-patterns that are
symptomatic of the long-range repetitive lattice nature of the crystal.
Sometimes these pictures seem so crystal-like themselves that, at first sight,
they might lead one to think that they rather directly mark the atomic
positions.  In fact, however, they display the symmetry of another lattice
that is dual (or reciprocal) to the one underlying the crystal structure. (See
Boxes \ref{CandP} and \ref{diffrac} for more on this).

For almost 80 years, the point-like feature of the diffraction image seemed to
be the characterizing property of crystals; so much so that the three concepts
of lattice symmetry, crystal structure, and pure point diffraction were
considered as synonymous. Thus it was a minor crisis for the field of
crystallography when in 1982 certain materials were found \cite{Dany} with
diffraction patterns that were as point-like as those of crystals, but showed
other symmetries that are not commensurate with lattice symmetry! So, these
new substances, which were definitely not crystals in the classical sense,
were quickly dubbed {\em quasi-crystals}, and opened a new branch of
crystallography.  At the same time, they brought forth a surge of new
mathematics with which to model the new geometry involved.

It is to this mathematical side that we turn in this article. For beyond the
many physical questions raised by these new quasicrystals, there is a bundle
of mathematical questions. What do we mean by `order', in particular by
`aperiodic order', how do we detect or quantify it, what do we mean by
repetition of patterns, what are the underlying symmetry concepts involved,
how can one construct well-ordered aperiodic patterns?  Beyond this, as one
quickly realizes, is the general question of how the new class of
quasicrystals and their geometric models are to be placed between the perfect
world of ideal crystals and the random world of amorphous or stochastic
disorder or, in other words, how can we characterize the level of `disorder'
that we may have reached?

\bigskip
\centerline{
\fbox{
\begin{minipage}{0.9\textwidth}
\centerline{\epsfxsize=\textwidth\epsfbox{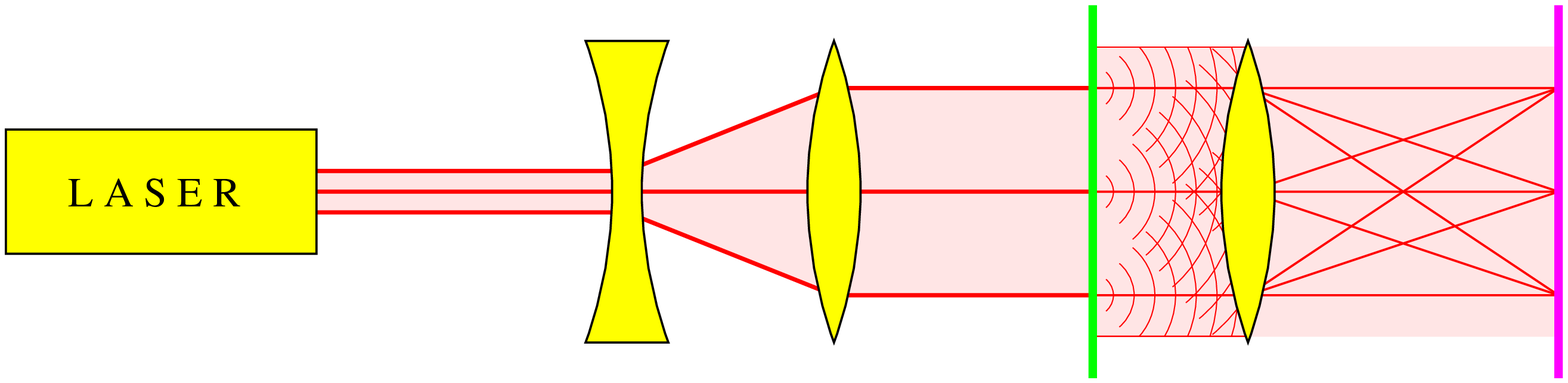}}
\begin{picbox}
  Experimental setup for optical diffraction\newline {\em The laser beam is
    widened by an arrangement of lenses and orthogonally illuminates the
    object located at the green plane.  The light that emanates from the
    object plane then interferes, and the diffraction pattern is given by the
    distribution of light that one would observe at an infinite distance from
    the object. By another lens, this pattern is mapped onto the pink plane.
    Whereas for a picture of the object, as for instance in a camera, light
    rays emanating from one point of the object ideally are focused again into
    a single point of the picture, the situation is different in diffraction
    --- light emanating from different regions within the object make up a
    single point of the diffraction pattern, as schematically indicated by the
    red lines in the right part of the figure.  Therefore the diffraction
    pattern carries information about the entire illuminated part of the
    object. It provides some kind of measure of the correlations, and thus an
    account of the degree of order, in the structure of the object.}
\label{laser}
\end{picbox}
\end{minipage}
}}
\bigskip

\section{Planar tilings}

A very instructive and also very attractive way to get a feeling for the ideas
involved is to look at two-dimensional tiling models. The two rhombi (the
so-called proto-tiles) shown in Box~\ref{pentiles} are clearly capable of
periodic stacking and so of lattice symmetry, the symmetry lattice being
generated by the two translational shifts shown. Another possibility is shown
below, which gives a tiling that is periodic in one direction and arbitrary
(in particular, possibly aperiodic) in the other. On the other hand, the
rhombi can also be used to tile the plane in the form of the famous Penrose
tiling, see Box~\ref{penfig}.

\bigskip
\centerline{
\fbox{
\begin{minipage}{0.9\textwidth}
\[
\centerline{\epsfxsize=0.5\textwidth\epsfbox{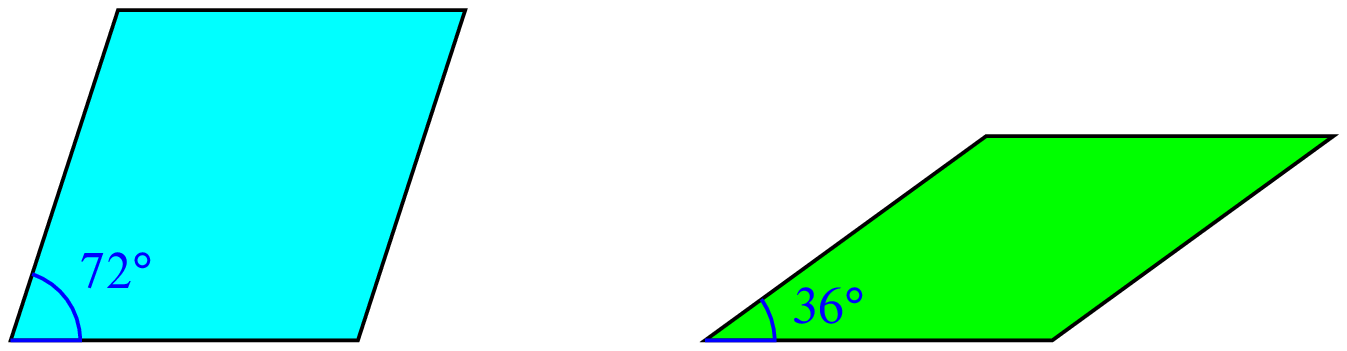}}
\]
\[
\centerline{\epsfxsize=0.95\textwidth\epsfbox{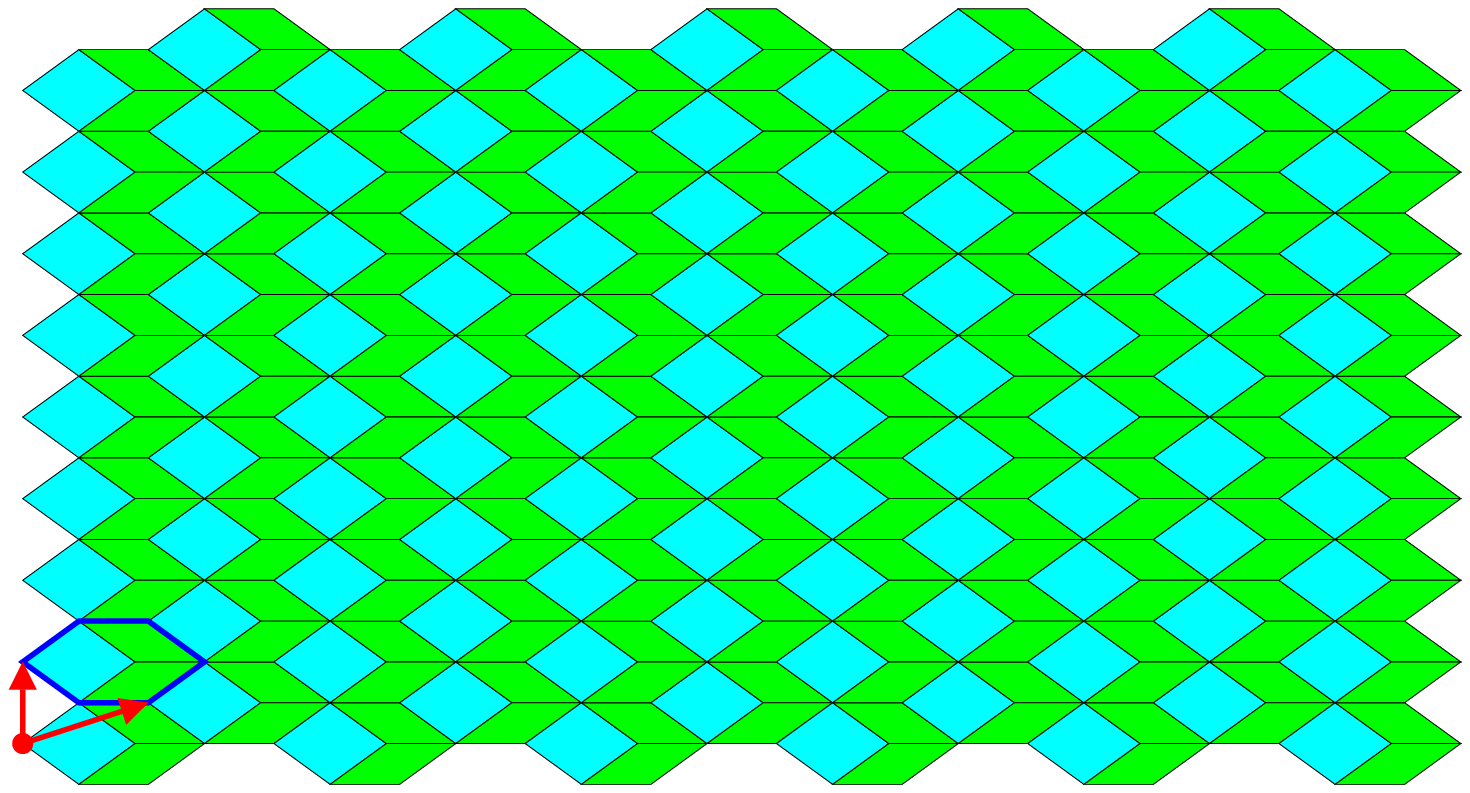}}
\]
\[
\centerline{\epsfxsize=0.95\textwidth\epsfbox{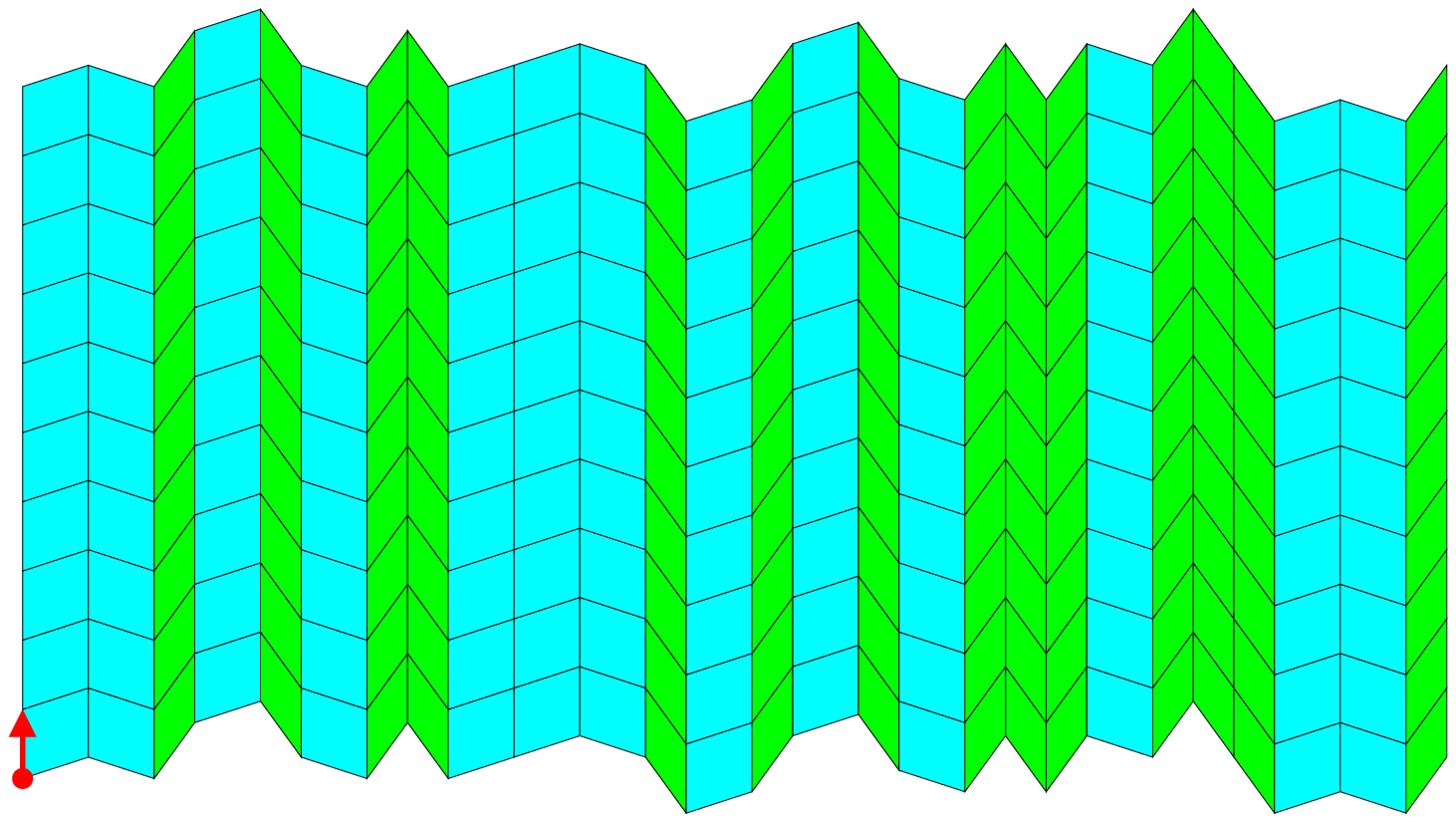}}
\]
\begin{picbox}
  The undecorated Penrose tiles and some of their assemblies\newline {\em The
    prototiles are two rhombi, a fat one with opening angle $72^\circ$ and a
    skinny one with $36^\circ$. They admit periodic arrangements like the one
    shown in the middle. The fundamental periods are indicated by arrows, and
    a fundamental domain in form of a hexagon is highlighted. It contains one
    fat and two skinny rhombi. Below, another arrangement is shown, which is
    periodic in the vertical direction, but admits an arbitrary `worm' of
    rhombi in the horizontal direction.  }
\label{pentiles}\smallskip
\end{picbox}
\end{minipage}
}}
\clearpage

\bigskip
\centerline{
\fbox{
\begin{minipage}{0.9\textwidth}
\[
\centerline{\epsfxsize=0.5\textwidth\epsfbox{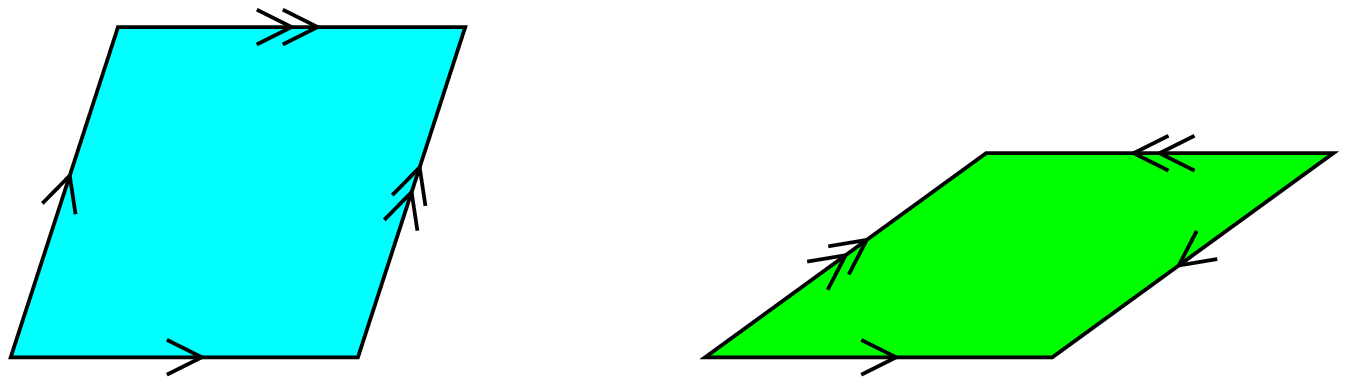}}
\]
\[
\centerline{\epsfxsize=\textwidth\epsfbox{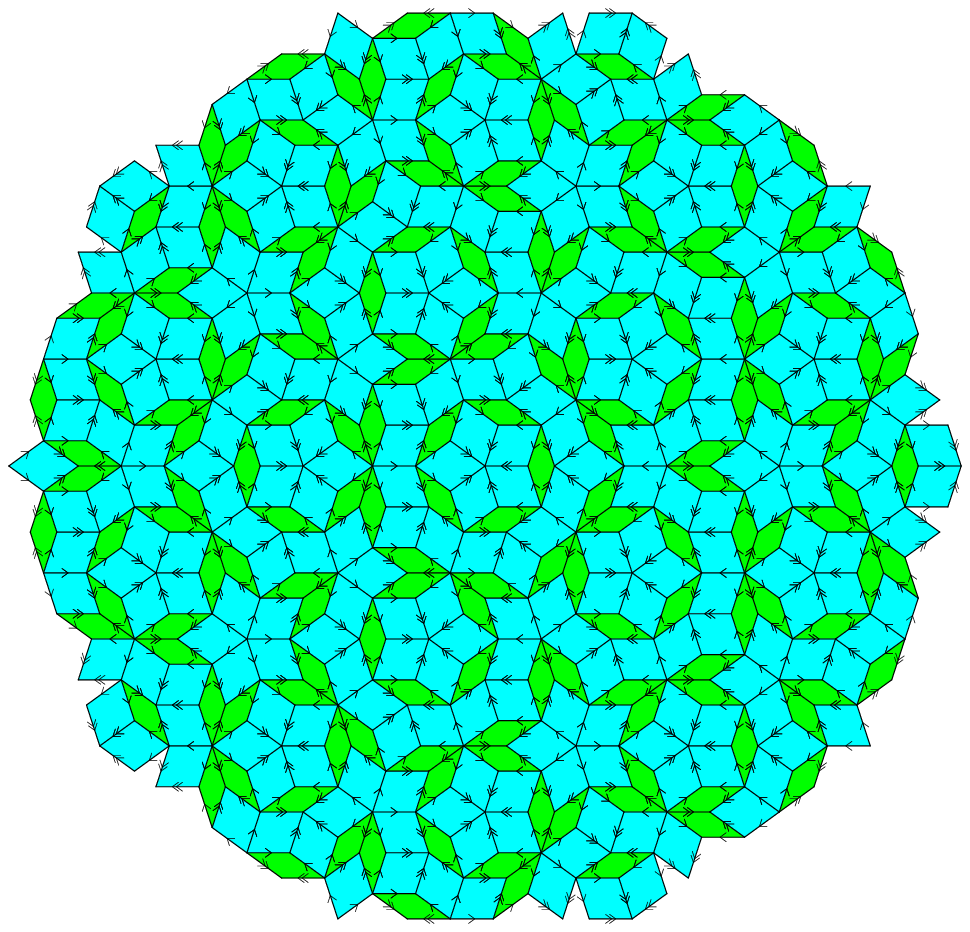}}
\]
\begin{picbox}
  A central patch of Penrose's aperiodic tiling\newline {\em The two rhombi of
    Box~\ref{pentiles} received a decoration of their edges by single and
    double arrows. If one now requires a perfect matching of all arrows on
    adjacent edges, the possible arrangements are highly restricted. In fact,
    the only permissible tilings of the entire plane are the so-called Penrose
    tilings. The different (global) possibilities cannot be distinguished by
    any local inspection.  A fivefold symmetric patch of such a tiling is
    shown above.  }
\label{penfig}\smallskip
\end{picbox}
\end{minipage}
}}
\clearpage

Part of the intriguing nature of the Penrose tiling, of which just a circular
fragment is shown in Box~\ref{penfig}, is the obvious question of what exactly
the rules might be for assembling these tiles. A properly constructed Penrose
tiling has several marvellous properties of which the two most important at
this point are:
\begin{itemize}
\item A complete Penrose tiling of the plane is
      strictly {\em aperiodic}  (in the  sense of being totally without 
      translational symmetries). Our particular example shows striking 
      five-fold symmetry.  
\item If we ignore the tiles and just look at their vertices instead (we might 
      think of the resulting point set as a toy model of an atomic layer) 
      then, remarkably, this set of points is itself pure point diffractive,
      i.e.\ in the optical bench of Box \ref{laser}, it produces a diffraction
      image on the screen with sharp spots only.
\end{itemize}

In Box~\ref{abpatch}, we see another aperiodic tiling, this time made out of
two very simple tile types, a square (which we actually dissect into two
isosceles triangles) and a rhombus. Its set of vertex points shows the same
type of diffraction image as the Penrose tiling, namely sharp spots only, this
time with eightfold symmetry (Box \ref{diffpatt}).  In Box \ref{abinf}, we see
the beautiful idea that is the secret behind many of the most interesting
tilings (including the Penrose tiles): the idea of inflating and subdividing.
To apply the idea here, we directly work with triangle and rhombus.

\bigskip
\centerline{
\fbox{
\begin{minipage}{0.9\textwidth}
\[
\centerline{\epsfxsize=0.5\textwidth\epsfbox{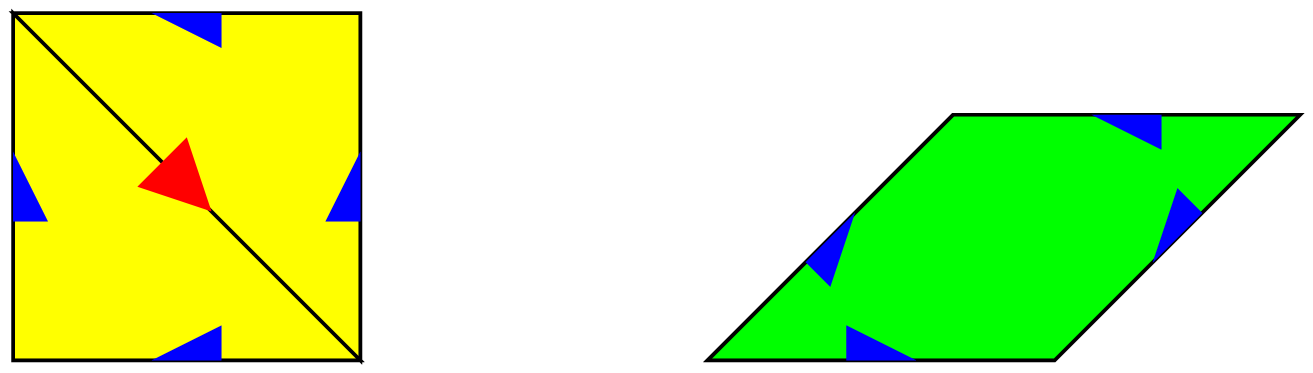}}
\]
\[
\centerline{\epsfxsize=\textwidth\epsfbox{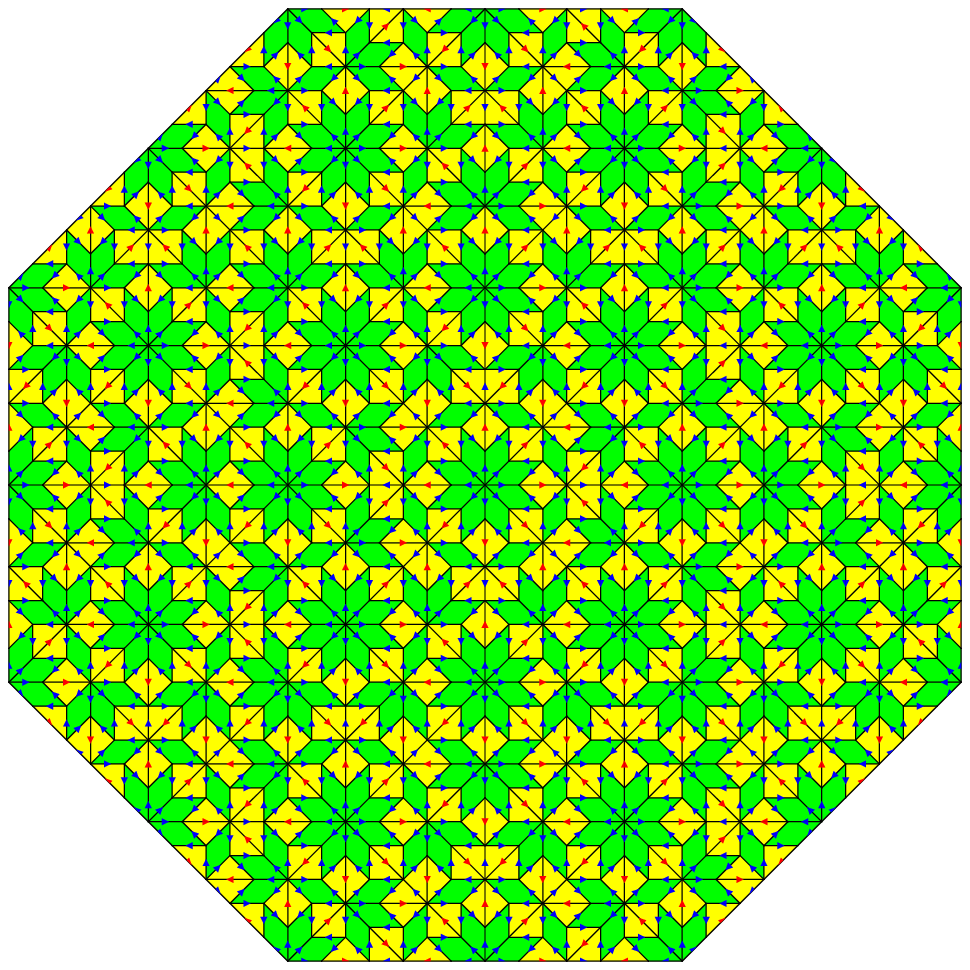}}
\]
\begin{picbox}
  A central patch of the octagonal Ammann-Beenker tiling\newline {\em The
    original prototiles are a square and a $45^\circ$ rhombus, decorated with
    blue arrows on the edges. For later use, the square is cut into two
    congruent isosceles triangles, carrying a red arrow on their common base.
    The orientation of arrows within each triangle is circular. Unlike the
    situation in the Penrose tiling, even with these arrows periodic tilings
    are still possible, for instance by repeating the square periodically. The
    octagonal patch shown belongs to the eightfold symmetric relative of the
    Penrose tiling, which is non-periodic and usually called the octagonal or
    the Ammann-Beenker tiling.  }
\label{abpatch}\smallskip
\end{picbox}
\end{minipage}
}}
\bigskip

\bigskip
\centerline{
\fbox{
\begin{minipage}{0.9\textwidth}
\[
\centerline{\epsfxsize=\textwidth\epsfbox{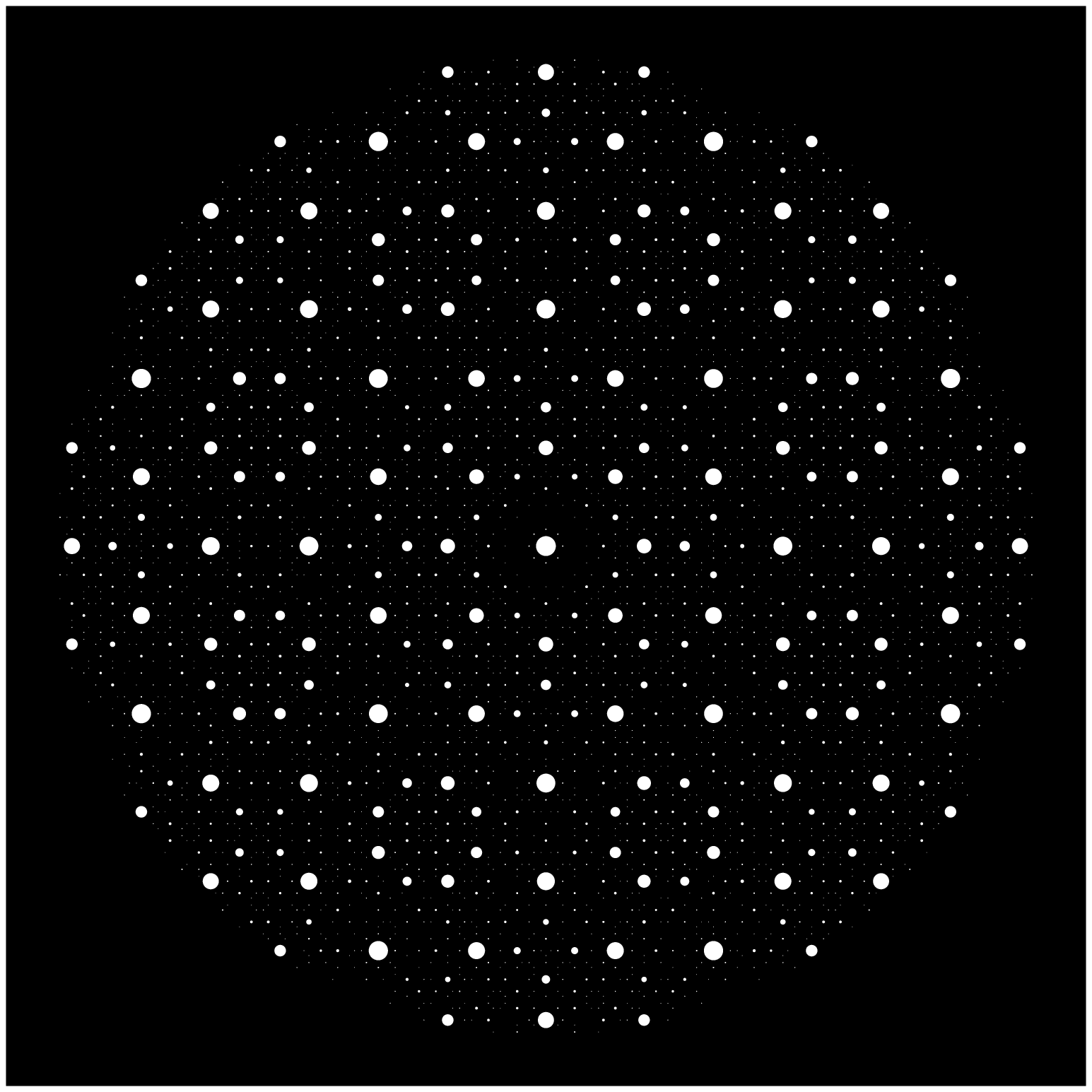}}
\]
\begin{picbox}
  Diffraction pattern\newline {\em Diffraction pattern of the octagonal
    Ammann-Beenker tiling.  The diffraction spots are indicated by circles
    whose area is proportional to the intensity of the diffraction peak.
    Spots with an intensity of less than 0.05\% of the intensity of the
    central spot have been discarded.  }
\label{diffpatt}\smallskip
\end{picbox}
\end{minipage}
}}
\bigskip

\bigskip
\centerline{
\fbox{
\begin{minipage}{0.9\textwidth}
\[
\centerline{\epsfxsize=\textwidth\epsfbox{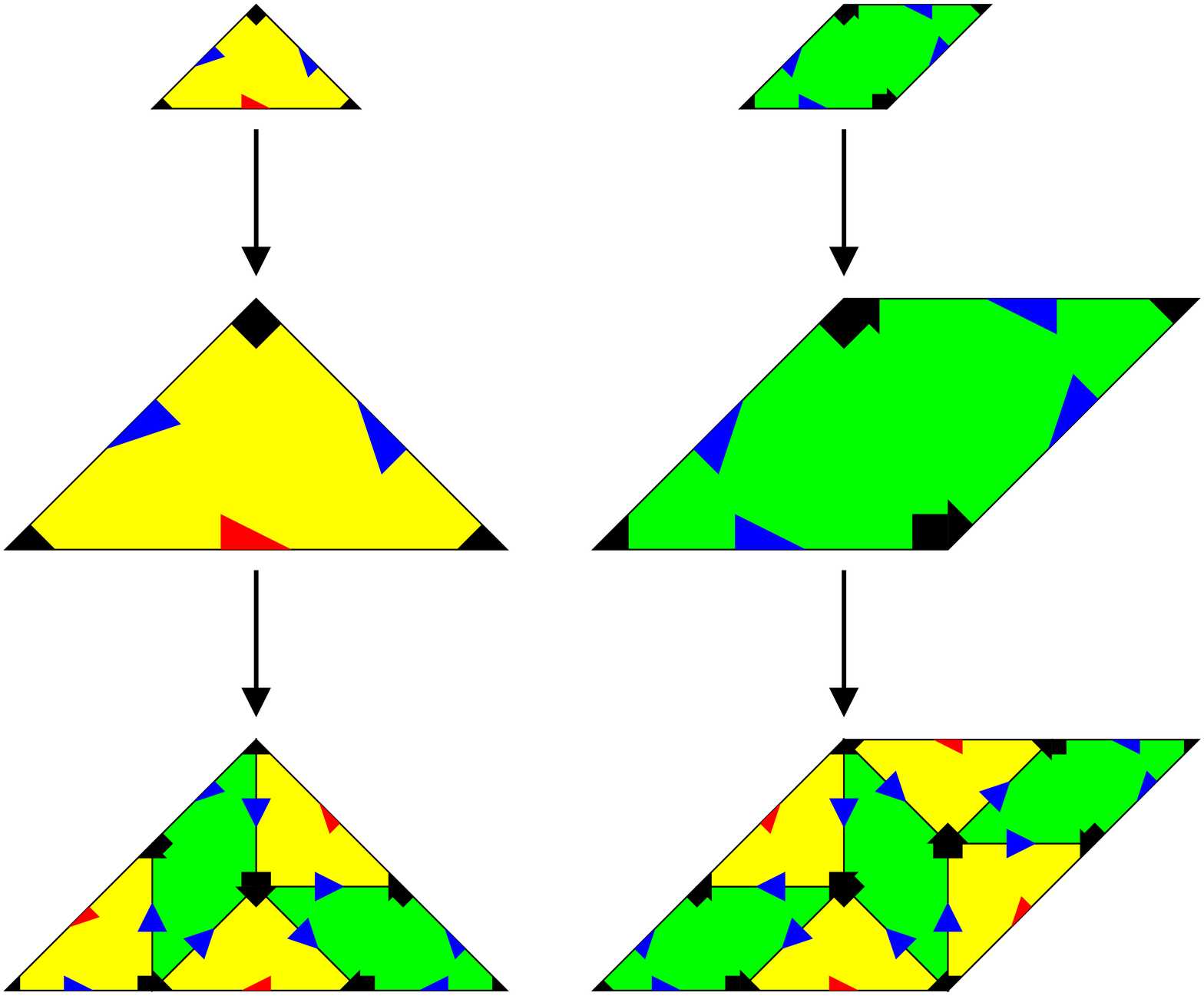}}
\]
\begin{picbox}
  Inflation rule for the octagonal Ammann-Beenker tiling\newline {\em The
    inflation procedure consists of two steps, a rescaling by a factor of
    $\alpha=1+\sqrt{2}$, followed by a dissection into tiles of the original
    size. In comparison to Box~\ref{abpatch}, corner markings have been added
    which break the reflection symmetry of the rhombus.  The patch shown in
    Box~\ref{abpatch} can be obtained by applying this inflation rule
    (ignoring the corner markings) to an initial patch that coincides with the
    central octagon, filled by eight squares and sixteen rhombi. The corner
    markings are vital for obtaining matching rules. A sequence of inflation
    steps starting from a single square is shown in Box~\ref{sequence}. Unlike
    the edge markings, and hence unlike the situation of the Penrose tiling,
    the corner markings cannot be reconstructed by local inspection of the
    undecorated tiling.  }\label{abinf}\smallskip
\end{picbox}
\end{minipage}
}}
\clearpage

The inflation scheme in Box~\ref{abinf} shows us how to inflate each tile by a
factor of $\alpha = 1 + \sqrt{2}$ and then how to decompose the resulting tile
into triangles and rhombi of the original size. With this new device, we have
a way of filling the whole plane with tiles. In comparison to
Box~\ref{abpatch}, we added some markers in the corners of the tiles which
will play some magic tricks for us later.  Starting from a single tile, or
from the combination of two triangles, and inflating repeatedly, we build up
the sequence as shown in Box~\ref{sequence}. Since there is no need to stop,
we may go on and do this forever.

It is now easy to see that the resulting octagonal tiling has an amazing
property: whatever finite pattern of tiles we see, that same pattern will be
repeated infinitely often, in fact we can even specify the maximum distance we
will have to search to find it again!  A pattern with such a property is
called repetitive. A perfect crystal is an example of a repetitive structure,
of course, but the inflation procedure produces interesting new cases.

How does this happen? Imagine the partial tiling obtained after $n$ inflations
of an original patch $P$ that consists of two triangles which build a square.
It is composed of triangle pairs and rhombi.  If we choose from it a patch
$P'$ which is a copy of $P$, then $n$ steps after this patch was created,
another patch $P''$ will show up which is a copy of $P'$. Furthermore, the
position and orientation of $P''$ relative to $P'$ will be the same as that of
$P'$ relative to the original $P$. Thus the pattern $P$, or a similar copy
thereof, is bound to appear over and over again. In our example, $P$ is just
made of two tiles, but this idea works for any patch $P$ that occurs somewhere
in the inflation process, no matter how big it is.

The reason behind this is that the square, centred at the origin, is the seed
of a fixed point under even numbers of inflation, as can be seen from the
sequence in Box~\ref{sequence}. The term `fixed point' means that the sequence
tends towards a global covering of the plane which is then left invariant
(hence fixed) by further pairwise inflation steps, i.e., we have reached a
stable pattern this way.

\bigskip
\centerline{
\fbox{
\begin{minipage}{0.9\textwidth}
\[
\centerline{\epsfxsize=\textwidth\epsfbox{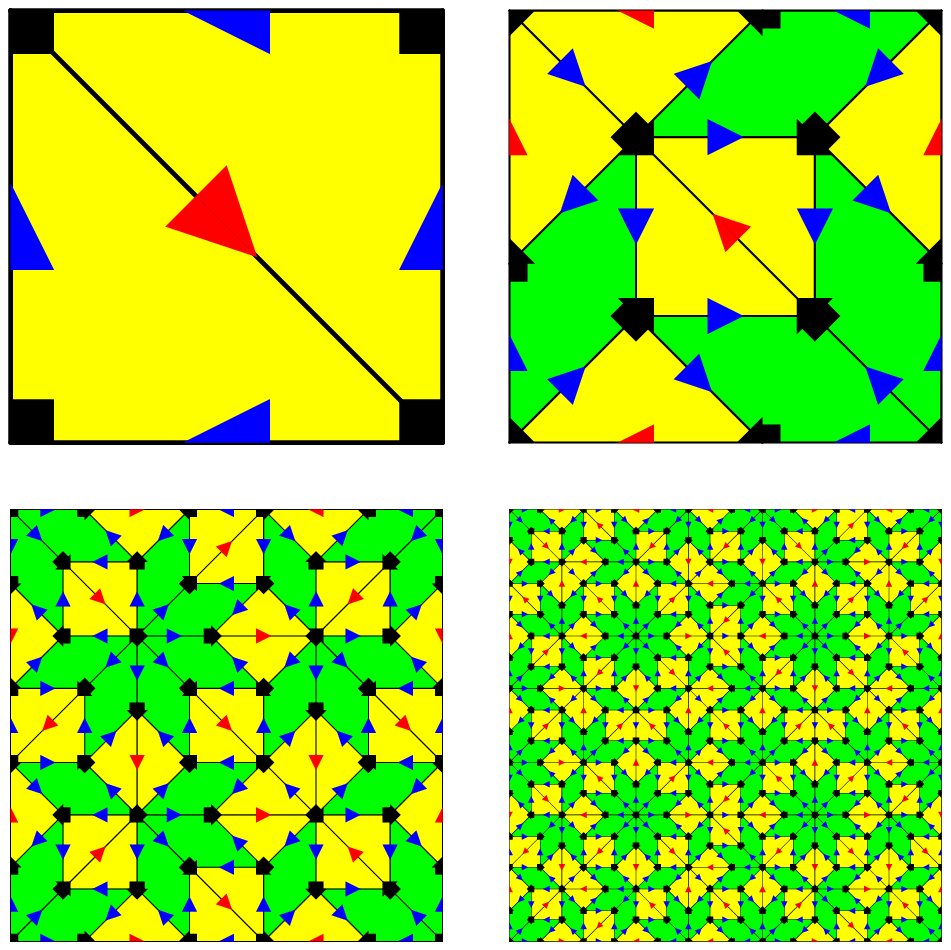}}
\]
\begin{picbox}
  Repeated inflation steps of the octagonal tiling\newline {\em The sequence
    shows a square as an initial patch and three successive applications of
    the inflation rule of Box~\ref{abinf}. (For the sake of presentability, we
    ignored the proper relative scale.) The inflation rule ensures that the
    corner markings always assemble a complete `house'. Alternatively,
    assembling patches tile by tile, all complete tilings of the plane with
    this property and matching arrows on all edges are locally
    indistinguishable from the fixed point tiling created by inflation. Thus,
    arrows and houses together establish perfect matching rules.  }
\label{sequence}\smallskip
\end{picbox}
\end{minipage}
}}
\clearpage

So our pattern is {\em repetitive}, but in fact it has no periodic component
at all! This is not self-evident yet, but it will become more so later. The
main point right now is that the tiling has the strange and seemingly
paradoxical property of having repetitivity on all scales, no matter how
large, but with no periodic repetition.  All patches repeat, but not
periodically!

The Penrose tilings can also be built through substitution and likewise are
repetitive without periodic repetition, see \cite{GS}.  Thus they too have the
striking property that you cannot really know where you are in the tiling by
looking at any finite region around you. It follows that it is not possible to
build such a tiling by any finite set of rules which tell you what to do next
by looking at some finite neighbourhood of your position!  To see why, imagine
that this were possible. Then every time the same pattern appeared, the rules
for continuing from it would be the same as those used for building at its
previous occurrence. The result is that the pattern would globally repeat.

Having said this, the next reaction is going to be that our next assertion
says the opposite.  In fact there are assignments of marks --- so-called
matching rules --- to the edges of the Penrose rhombi (Box~\ref{penfig}), or
to the edges and corners of the Ammann-Beenker tiles (Boxes~\ref{abinf} and
\ref{sequence}), such that, if they are match everywhere in the tiling, the
result is a perfect Penrose or a perfect Ammann-Beenker tiling, respectively.
What is the catch?

The problem is that these matching rules guarantee that what you are getting
is a Penrose tiling {\em as long as you never get stuck}. The trouble is that
to not get stuck requires knowledge of the entire tiling to that point --- it
is not derivable from local information only!

\section{Cut and project sets}
 
In view of these difficulties, one might ask what other possibilities exist to
systematically create arbitrarily large faultless patches of these tilings.
The idea of what is going on is more easily understood by first considering an
even simpler object, namely a one-dimensional inflation tiling.  This time we
begin with two tiles
\[
\centerline{\epsfbox{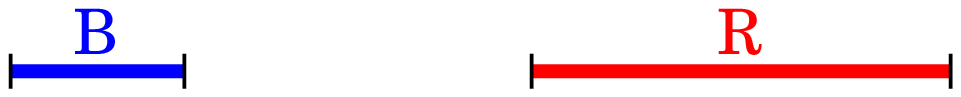}}
\]
which we call B (for blue) and R (for red), respectively. We give the short
tile B the length $1$ and the long tile R the length $\alpha=1+\sqrt{2}$ (the
same number also appears in the octagonal tiling). Inflation is stretching by
a factor of $\alpha$, followed by a subdivision which is consistent with
$\alpha\cdot 1=\alpha$ and $\alpha\cdot\alpha=2\alpha+1$. The final result is
\[
\centerline{\epsfxsize=\textwidth\epsfbox{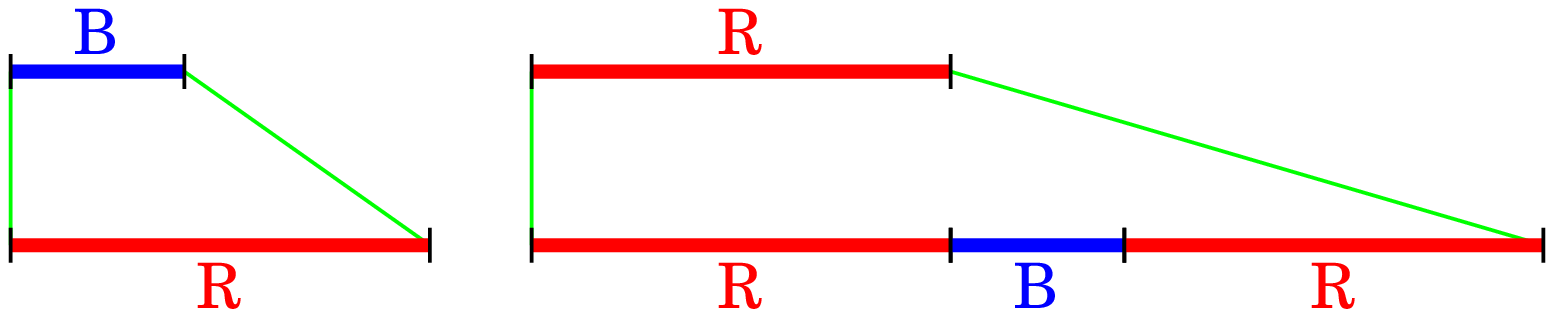}}
\]
Starting from a pair of R-tiles, centred at the origin, we have successively
\[
\centerline{\epsfxsize=\textwidth\epsfbox{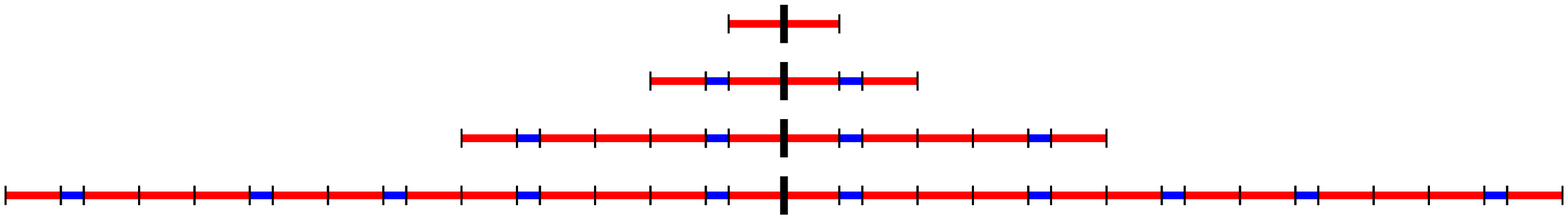}}
\]
Using coordinates to label the left end point of each
tile we have
\[
\centerline{\epsfxsize=\textwidth\epsfbox{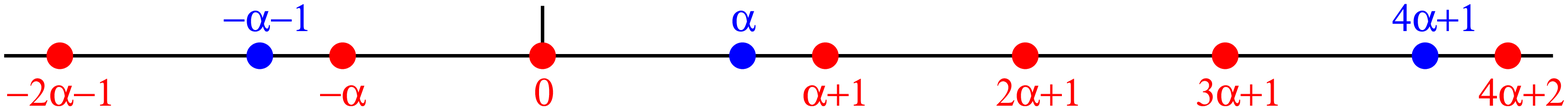}}
\]
The corresponding points form an infinite set $A = \{ \dots -\alpha -1,
-\alpha, 0 , \alpha, \alpha+1, 2\alpha +1, ...\}$.

What is striking about the points of $A$ is that they are all of the form $u +
v\sqrt 2$. How can we see which points $u+v\sqrt{2}$ are present and which
not? Everyone knows that it is a good idea in dealing with square roots to see
what happens if you change the sign of the square root. (Think of the high
school exercises in rationalizing expressions of the form $\frac{1}{1+\sqrt
  2}$.)

Let us use this trick of replacing each appearance of $\sqrt 2$ by its
conjugate, $-\sqrt{2}$. This conjugation is called the star map, the image of
a point $x=u+v\sqrt{2}$ is $x^{*}=u-v\sqrt{2}$.  Box~\ref{CandP} shows a plot
of our points. We make a new picture in which each point $x$ is ``lifted'' to
the point $(x, x^*)$ in the plane. Our points of interest are shown against a
backdrop consisting of all possible points $(u+v\sqrt 2,u-v\sqrt 2)$ where
$u,v$ are integers.

\bigskip
\centerline{
\fbox{
\begin{minipage}{0.9\textwidth}
\[
\centerline{\epsfxsize=\textwidth\epsfbox{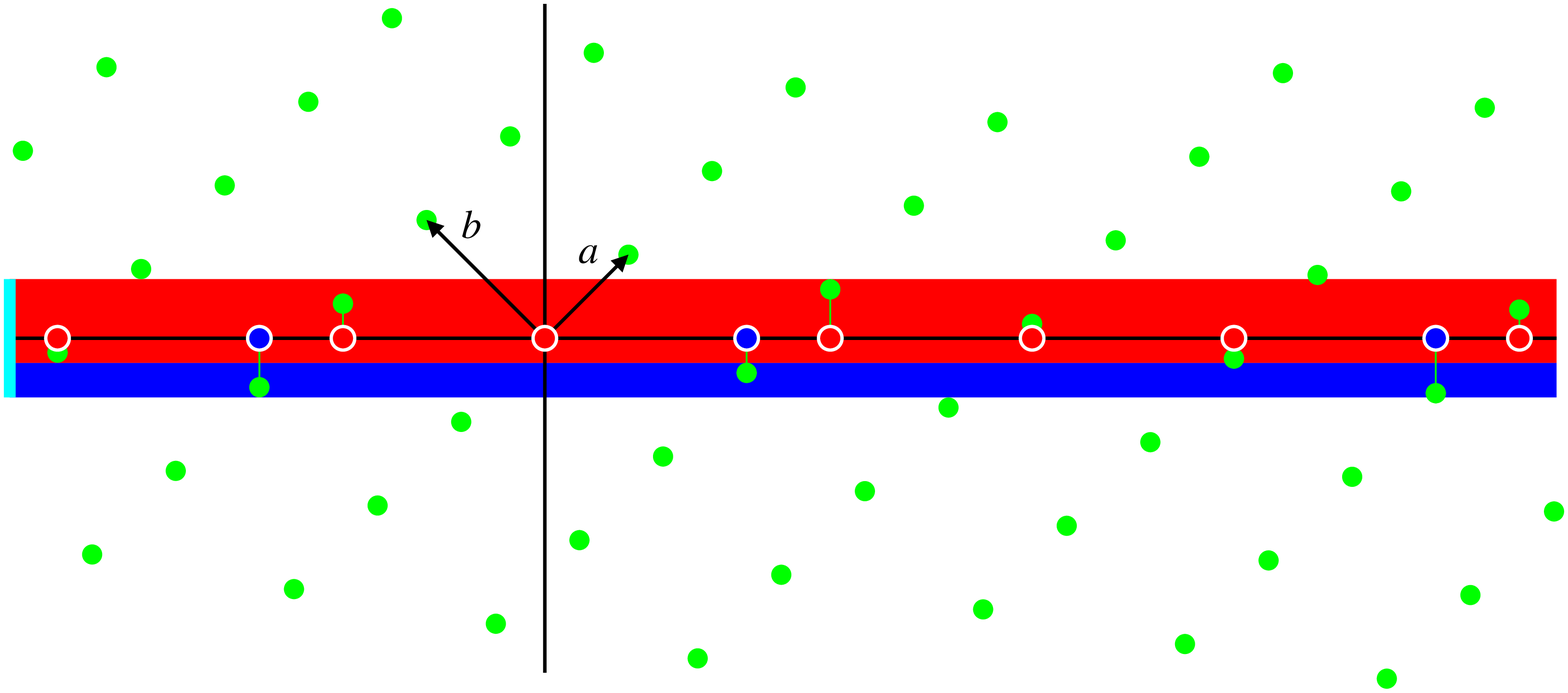}}
\]
\begin{picbox}
  An alternative way to construct the point set $A$\newline {\em The green
    points form the lattice $\{(u+v\sqrt{2},u-v\sqrt{2})\mid \mbox{$u,v$
      integer}\}$ which is spanned by the basis vectors $a$ and $b$.  The
    orientation of the strip is irrational with respect to the lattice, i.e.,
    the black line at its centre hits the origin, but no further lattice
    point. The green points within the strip are orthogonally projected onto
    the horizontal black line and are coloured according to their vertical
    position in the strip. The resulting set of red and blue points coincides
    with the point set constructed above by inflation.  }
\label{CandP}\smallskip
\end{picbox}
\end{minipage}
}}
\bigskip\bigskip\bigskip

The effect is striking. The entire set of points, including the backdrop,
produces a lattice (a mathematical crystal). The B and R points now appear in
a band that runs from height $-\frac{1}{\sqrt 2}$ to $\frac{1}{\sqrt 2}$.
Furthermore, the B points come from the bottom portion of the band, from
$-\frac{1}{\sqrt 2}$ to $\frac{1}{\sqrt 2} -1$, and the R points from the
remaining top portion of the band. The actual points labelling our tiling,
i.e.\ the set $A$, can be obtained just by dropping the second coordinate of
each lattice point that lies in the band --- in other words by projecting it
onto the horizontal axis.

Now one sees that it is incredibly easy to compute the left hand end points of
our $1$D tiling, and hence to get hold of the tiling itself. On a computer,
generate, in some ordered way, points of the type $u+ v\sqrt 2$. For each one
look at its conjugate $u- v\sqrt 2$.  Test whether this number lies in either
of the intervals corresponding to B and R points (e.g., $-\frac{1}{\sqrt{2}} <
u - v \sqrt{2} < \frac{1}{\sqrt{2}} $ for B points) and choose the point and
its colour accordingly. What we have accomplished here, apart from the visual
clarity, is a remarkable way of connecting the geometry of our tiling with an
algebraic method of calculating it.

A point set that can be described in this way (by cutting through a lattice
and projecting the selected points) is called, not surprisingly, a cut and
project set. In this case the object that is used to cut (or to sweep out) the
correct band is the vertical line segment indicated in black in
Box~\ref{CandP}. It is called the {\em window\/} of the projection method.

Another benefit of the cut and project view is that it shows immediately why
the resulting point sets are aperiodic. For example, a period of our set of
red and blue points is a shift $t$ (to the left or right) that moves the set
on top of itself. Necessarily it would be of the form $r + s\sqrt 2$ since all
our points have this form. However, after our lift into $2$-space, we would
then find that shifting by $(r + s\sqrt 2,r - s\sqrt 2)$ takes the strip onto
itself!  This is impossible unless $r - s\sqrt 2 =0$, i.e., $r = s \sqrt{2}$.
However, $\sqrt{2}$ is irrational, while $s,r$ are integers, so the only
solution is $r=s=0$, and the only period is $0$.

\section{The projection approach to planar tilings}

The octagonal tiling, or more precisely the positions of its vertices, can
also be described as a cut and project set. This goes via the projection of
the points of a certain lattice in four dimensions, swept out by an octagon.
We explain this in more detail.

The initial pool of points from which we select is given by the set $M$ of all
integer linear combinations
$\{u^{}_1a^{}_1+u^{}_2a^{}_2+u^{}_3a^{}_3+u^{}_4a^{}_4\mid
\mbox{$u^{}_1,u^{}_2,u^{}_3,u^{}_4$ integer}\}$ of the four unit vectors shown
in left diagram of Box~\ref{stars}. This is a dense point set in the plane,
and it is the two-dimensional analogue of the set $\{u+v\sqrt{2}\mid
\mbox{$u,v$ integer}\}$ used above. Since the octagonal tiling consists of
squares and rhombi (with unit edge length, say), the distance between any two
vertex points is of this form, i.e.\ an element of $M$. Also the star map has
an analogue, and it comes about simply by replacing the four vectors of the
left diagram by those of the right diagram of Box~\ref{stars}; that is,
$x=u^{}_1a^{}_1+u^{}_2a^{}_2+u^{}_3a^{}_3+u^{}_4a^{}_4$ is mapped to
$x^{*}=u^{}_1a_1^{*}+u^{}_2a_2^{*}+u^{}_3a_3^{*}+u^{}_4a_4^{*}$. As before,
the set of pairs $(x,x^{*})$ forms a lattice, this time in four dimensions.

\bigskip
\centerline{
\fbox{
\begin{minipage}{0.9\textwidth}
\[
\centerline{\epsfxsize=\textwidth\epsfbox{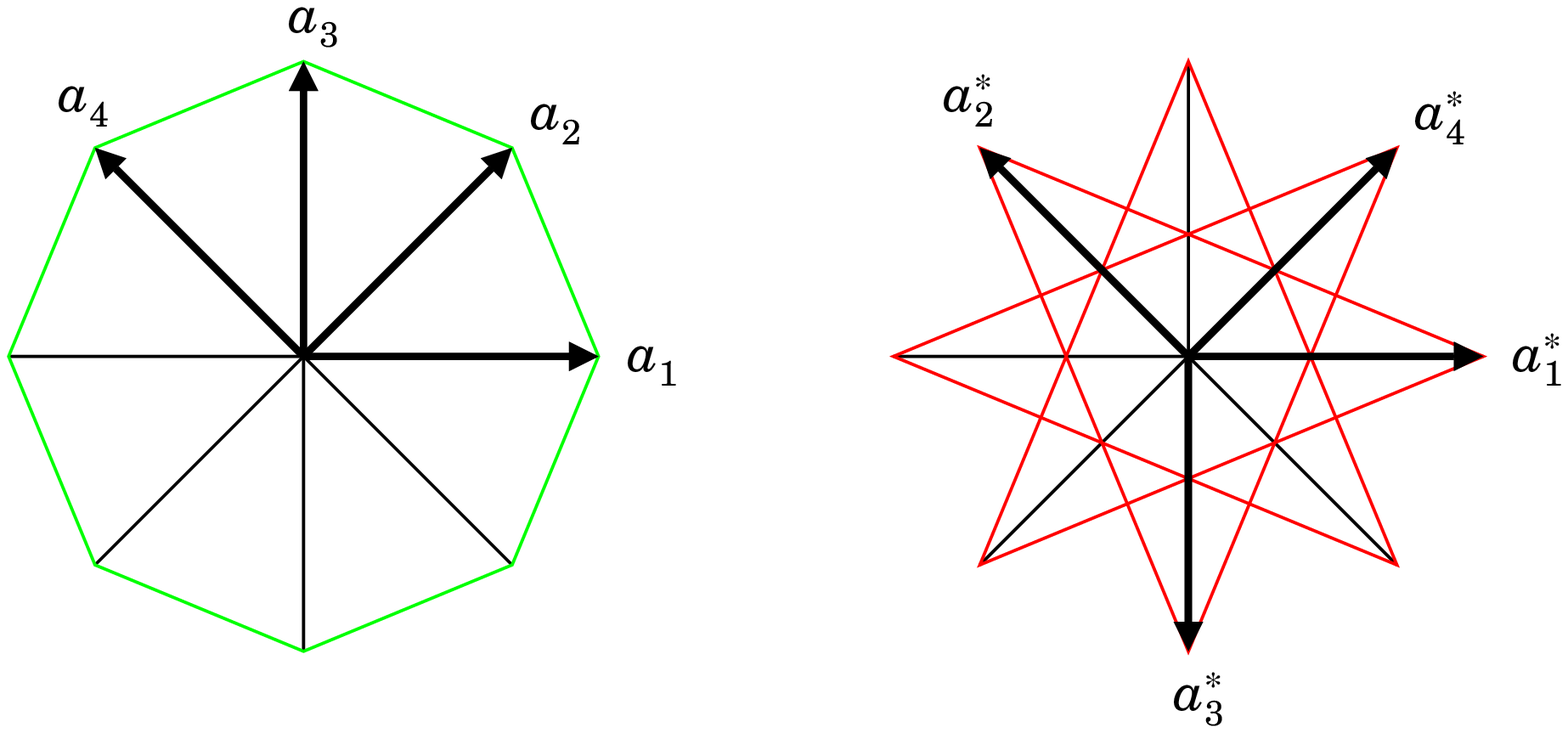}}
\]
\begin{picbox}
  The two ways to count to four (and hence to eight)\newline {\em The two sets
    of vectors used to construct the octagonal tiling, $a_{i}$ (left, for
    tiling space) and $a_{i}^{*}$ (right, for internal space), $i=1,2,3,4$.
    The change from $a_{i}$ to $a_{i}^{*}$ demonstrates the action of the
    $*$-map in this case.}
\label{stars}\smallskip
\end{picbox}
\end{minipage}
}}
\bigskip\bigskip\bigskip

The vertex set of the Ammann-Beenker tiling can now be given as the set of
points $x$ whose image $x^{*}$ under the star map lies inside a regular
octagon of unit edge length. We can now link this back to our previous
approach via inflation. If we start from a unit square and keep on inflating,
as shown in Box~\ref{sequence}, the images of the vertex points under the star
map will densely populate this octagon in a uniform way, see Box~\ref{weyl}.

Needless to say, the additional visual clarity obtained from a $4$D
description is debatable!  Still, the conceptual idea is very powerful,
providing the essential link between geometry, algebra, and analysis that is
at the heart of much of our understanding of aperiodic order.

Likewise the points of the Penrose tiling can be given a cut and project
interpretation, as do many other similar pointsets. In both cases, the
aperiodicity can be shown in the same way as for our one-dimensional example.

Another tiling of physical interest is built from a square and an equilateral
triangle. The example shown in Box~\ref{dodeca} can be created by a slightly
more complicated inflation rule, or alternatively once again by the cut and
project method. In this case, however, the corresponding window shows a new
feature: it is a compact set with fractal boundary. An approximation is also
shown in Box~\ref{dodeca}.

\bigskip
\centerline{
\fbox{
\begin{minipage}{0.9\textwidth}
\[
\centerline{\epsfxsize=0.5\textwidth\epsfbox{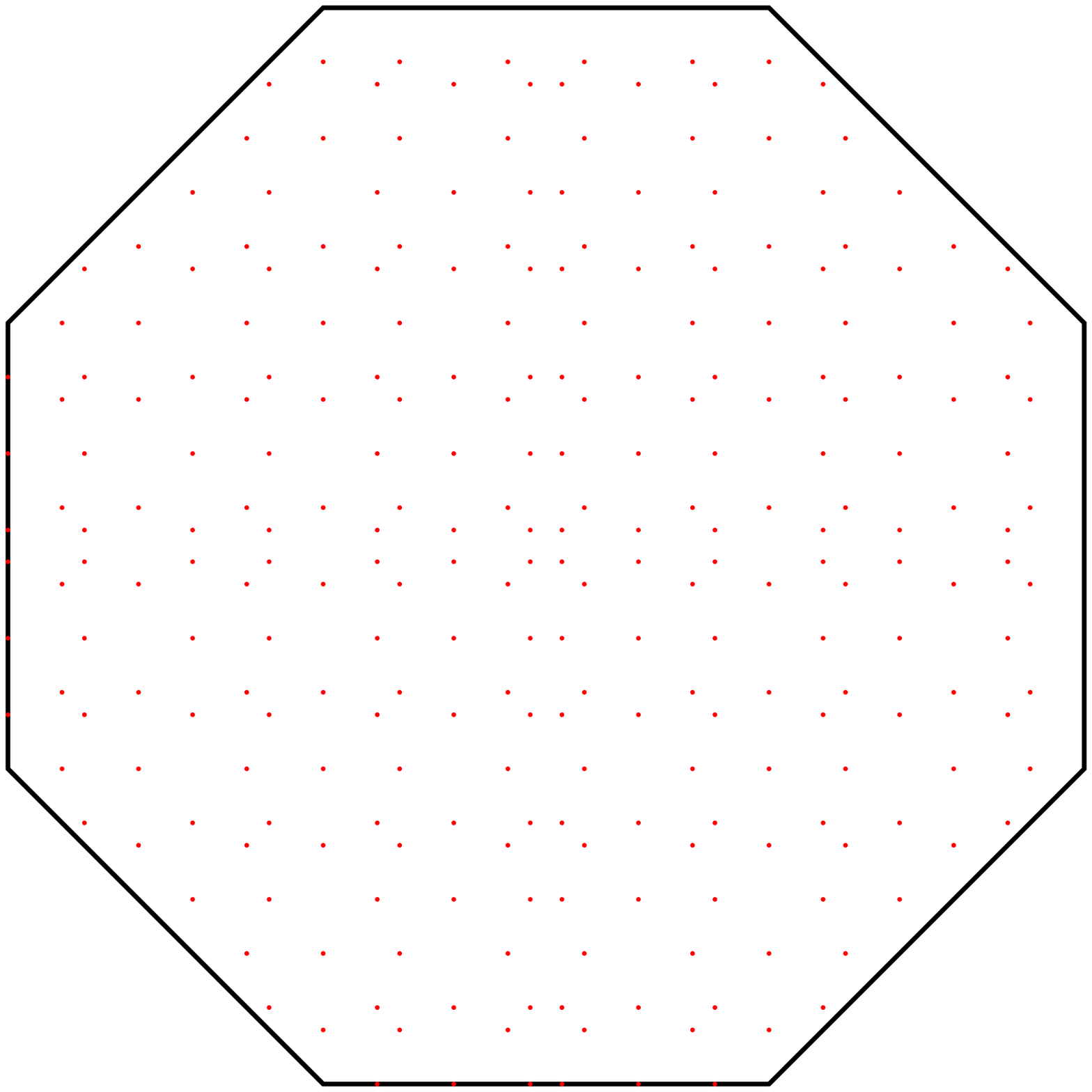}
\epsfxsize=0.5\textwidth\epsfbox{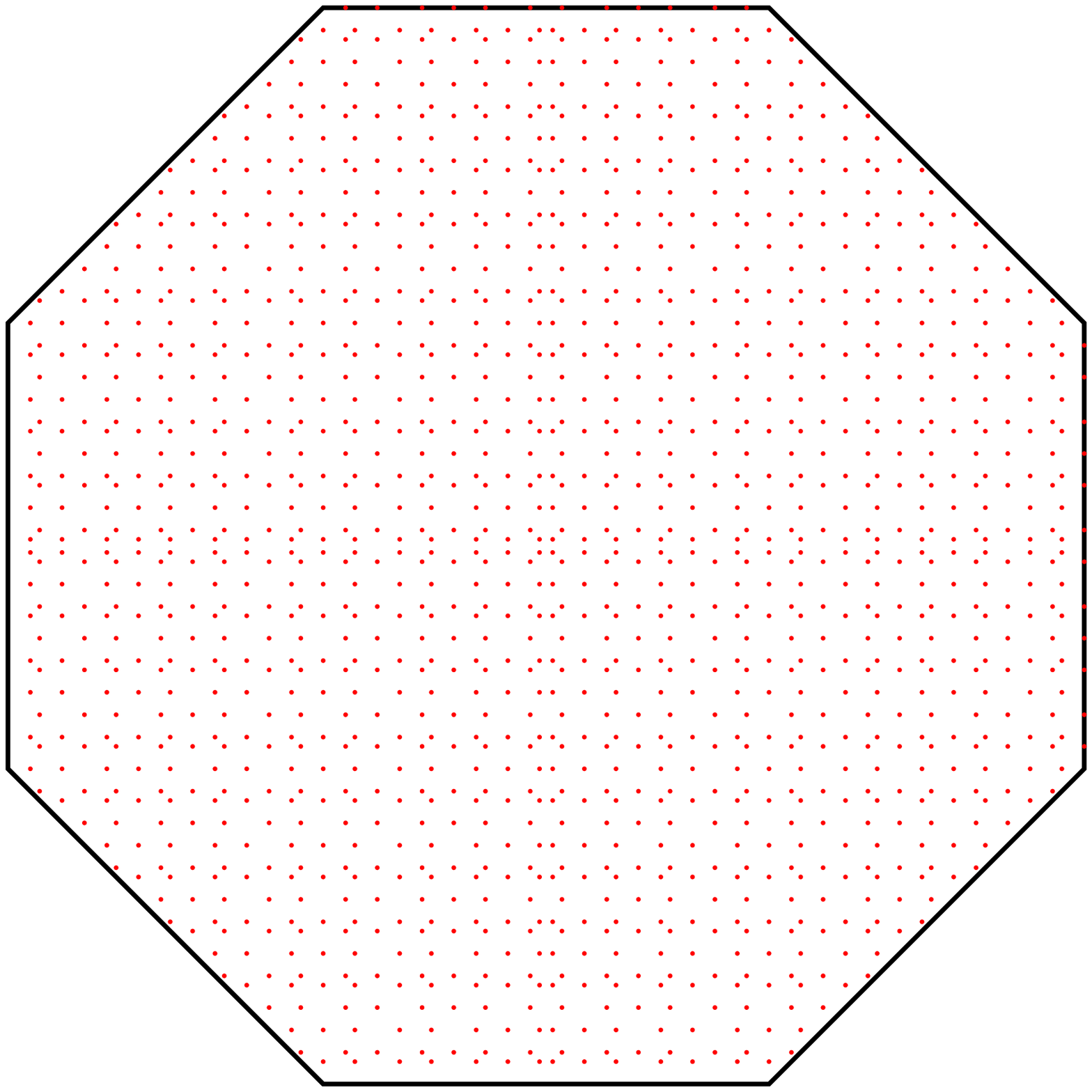}}
\]
\[
\centerline{\epsfxsize=0.5\textwidth\epsfbox{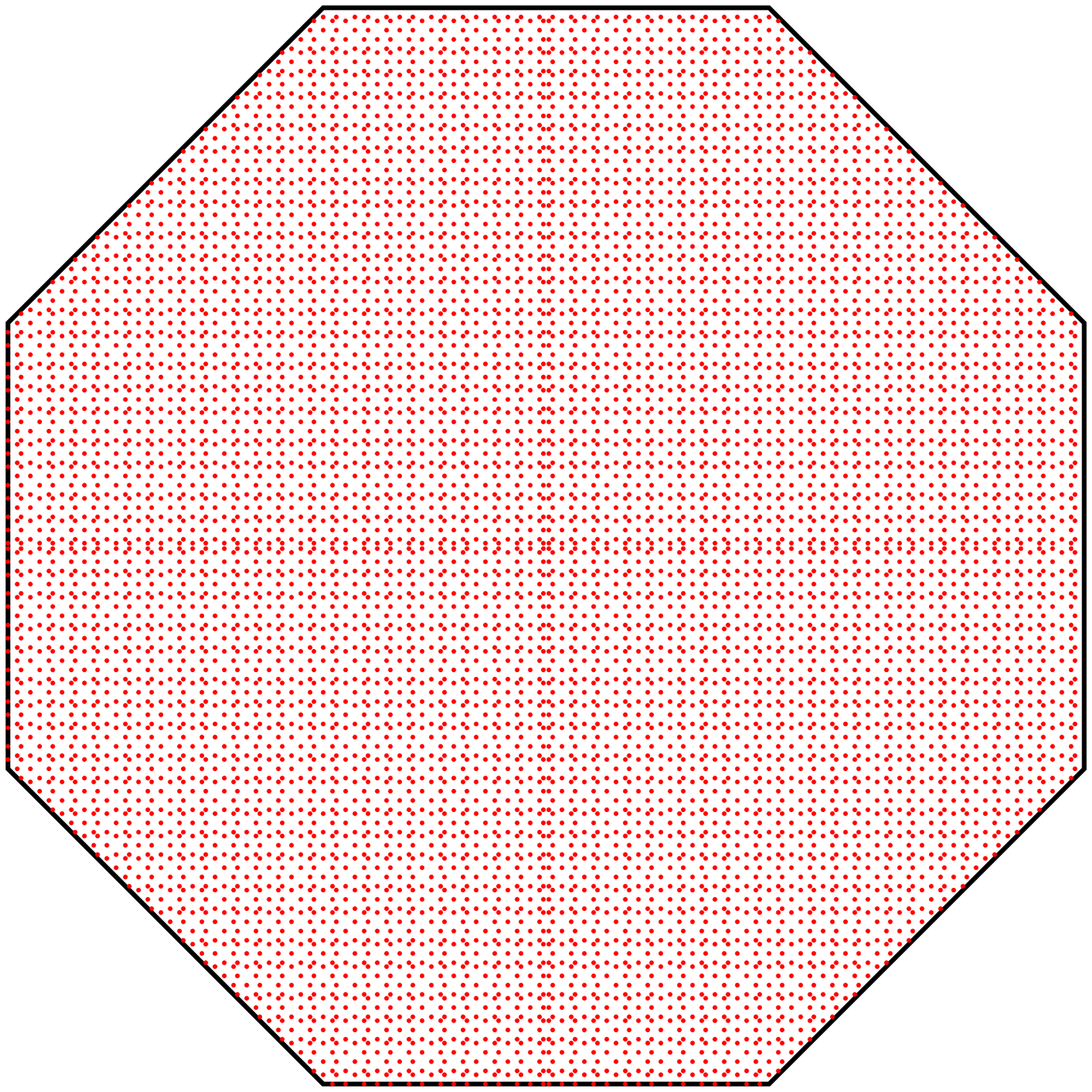}
\epsfxsize=0.5\textwidth\epsfbox{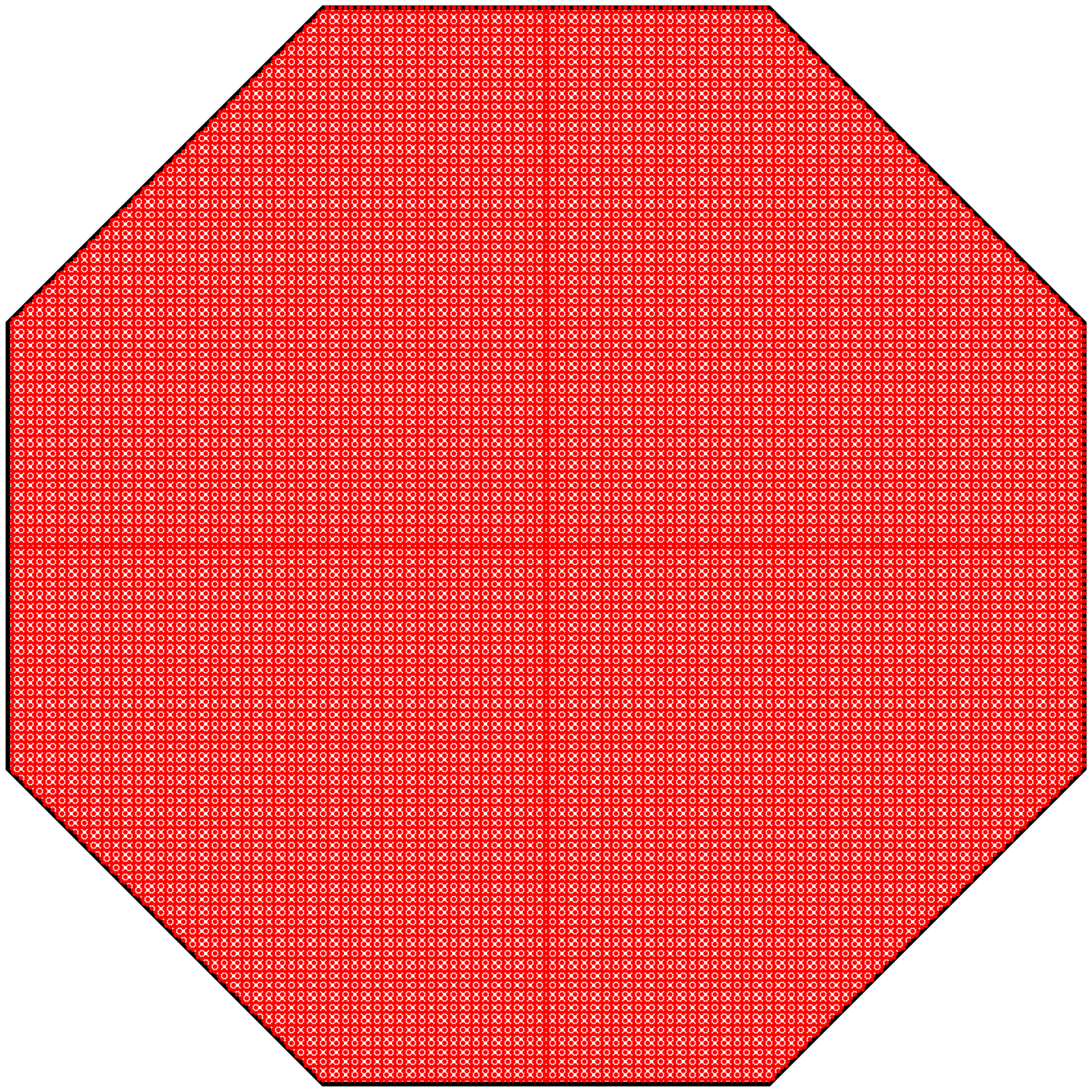}}
\]
\begin{picbox}
  Filling the octagon in internal space\newline {\em The image points $x^{*}$
    under the star map of the vertex points are shown for larger and larger
    patches of the octagonal tiling, obtained by inflation of a square as
    shown in Box \ref{sequence}. Eventually, the points populate the regular
    octagon with uniform density. Here, the first picture of the sequence
    corresponds to the largest patch of Box \ref{sequence}.}
\label{weyl}\smallskip
\end{picbox}
\end{minipage}
}}
\bigskip

\bigskip
\centerline{
\fbox{
\begin{minipage}{0.9\textwidth}
\[
\centerline{\epsfxsize=0.65\textwidth\epsfbox{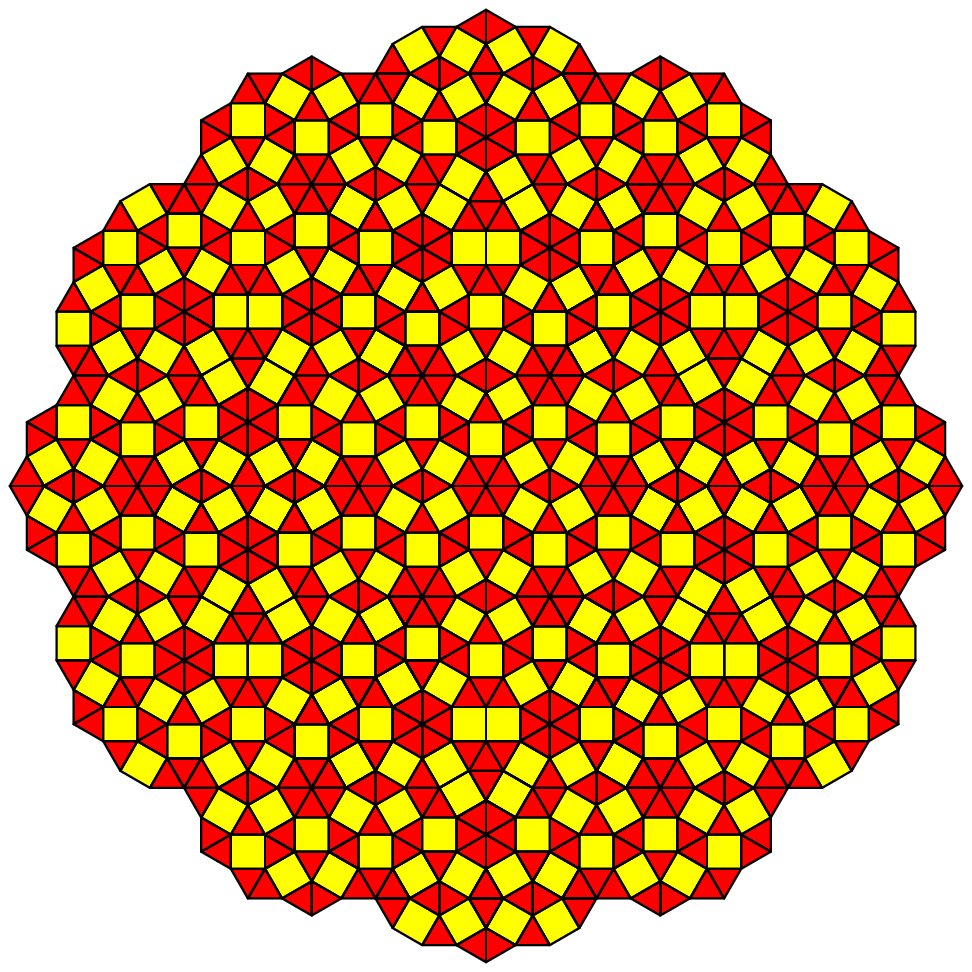}}
\]
\[
\centerline{\epsfxsize=0.65\textwidth\epsfbox{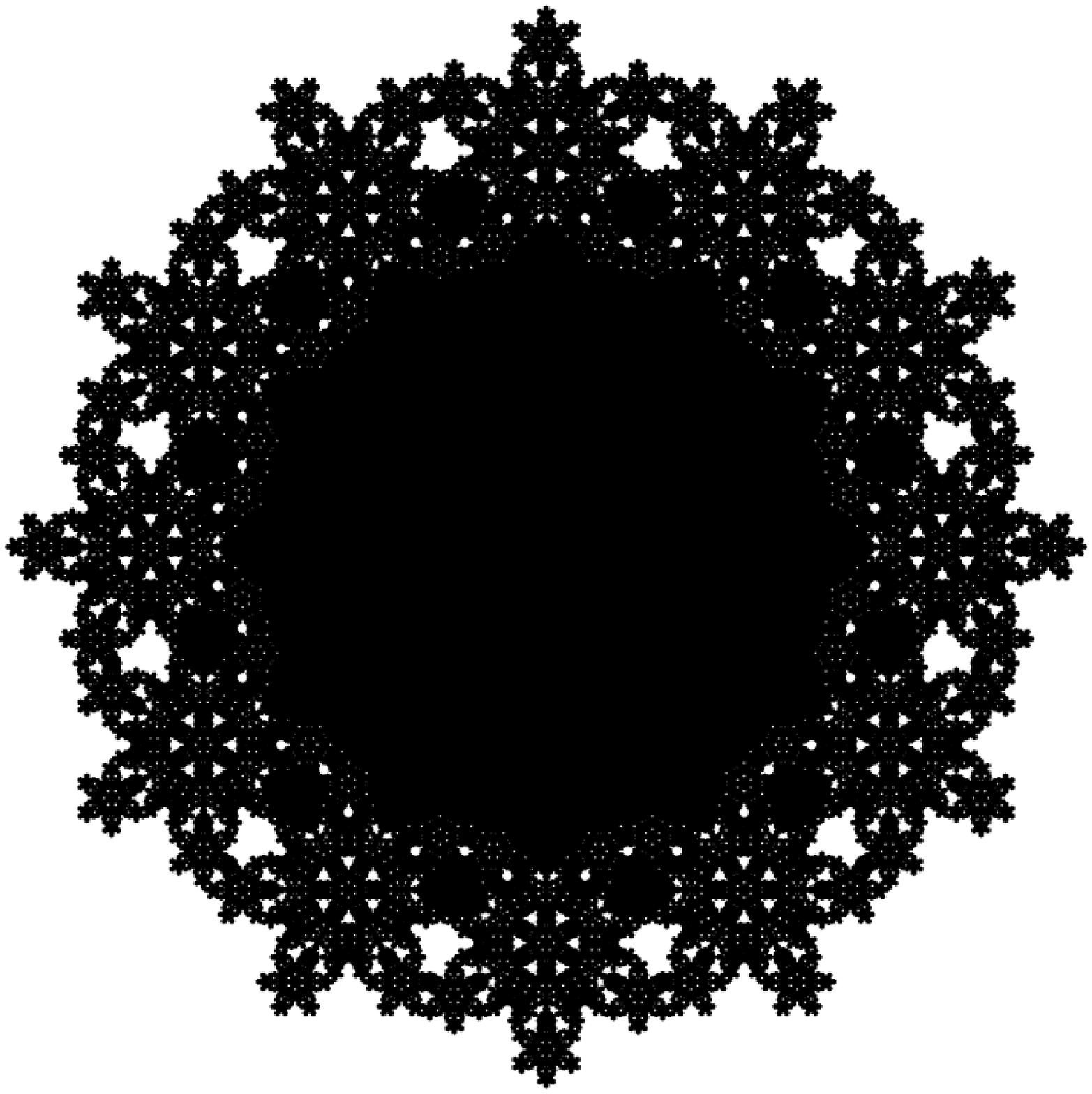}}
\]
\begin{picbox}
  Quasiperiodic square triangle tiling\newline {\em This example of a
    square-triangle tiling can either be obtained by an inflation rule or by
    projection from a lattice in four dimensions. The points selected for
    projection lie in a generalized `strip' whose cross section is a
    twelvefold symmetric object with fractal boundary.  }
\label{dodeca}\smallskip
\end{picbox}
\end{minipage}
}}
\bigskip

\section{The origin of diffraction}

The picture that we see in Box~\ref{CandP} offers us considerable insight into
the diffractive nature of sets that can be described as cut and project sets.
The background is a lattice (crystal) and this, from the classical theory of
crystals, is supposed to have perfect diffraction, i.e., the entire
diffraction image is composed of sharp peaks only. The trick is how to
restrict this down to the points in the band and ultimately to our line of
points. Box~\ref{diffrac} shows a picture of what happens. The bottom figure,
which looks like an irregular comb, shows the diffraction of the points $A$ of
our $1$D tiling.  The diffraction intensity is shown here not by the size of
the dots, but rather by the length of the teeth of the comb.

Above it is the diffraction picture of the background lattice, another
lattice, that, as we mentioned before, is called the dual lattice. The points
that carry the teeth of the comb (i.e. the spots of the diffraction) are
nothing other than the projections of the points of the dual lattice --- and
this time {\em all\/} of them.  The lengths of the teeth are provided by the
profile on the right hand side.  Where that profile comes from is a longer
story. (Engineers may recognize its similarity to the Fourier transform of a
single square pulse. It is, in fact, the square of the Fourier transform of
the characteristic function of the interval defining the band.)

The teeth of the comb lie actually dense on the line. However, due to the
damping nature of the profile, most of them are so small that, no matter what
finite resolution we may use, we can see only a small fraction of them, and
hence only an effectively discrete set of teeth, or spots, as in Box
\ref{diffpatt}.

\bigskip
\centerline{
\fbox{
\begin{minipage}{0.9\textwidth}
\[
\centerline{\epsfxsize=\textwidth\epsfbox{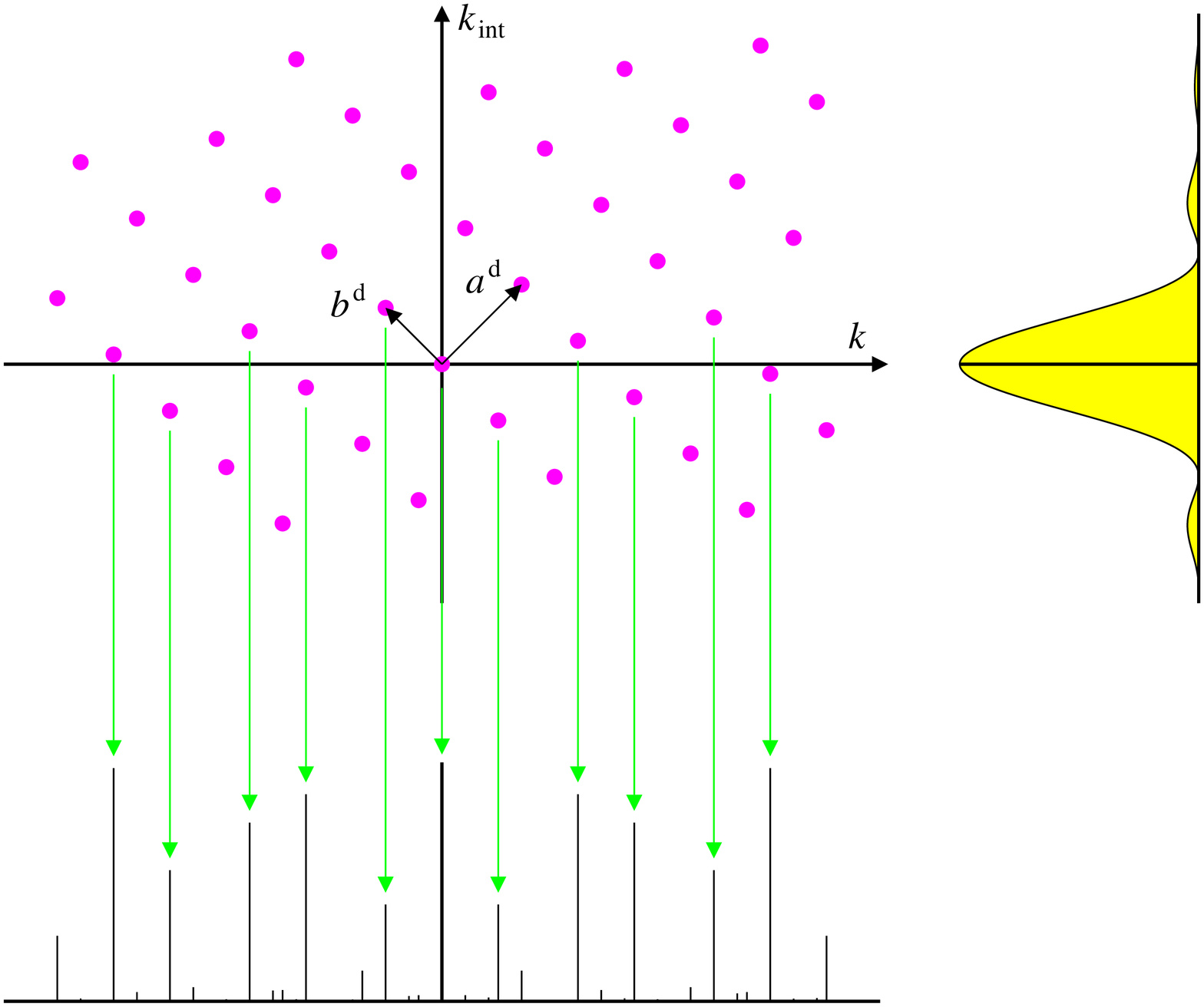}}
\]
\begin{picbox}
  Explanation of the diffraction pattern\newline {\em The pink points indicate
    the lattice dual to the lattice of Box~\ref{CandP}. It is explicitly given
    by
    $\{(\frac{m}{2}+\frac{n\sqrt{2}}{4},\frac{m}{2}-\frac{n\sqrt{2}}{4})\mid
    \mbox{$m,n$ integer}\}$. The lattice is spanned by the vectors $a^{\rm d}$
    and $b^{\rm d}$ which satisfy the scalar product relations $a^{\rm d}
    \cdot a = b^{\rm d} \cdot b = 1$ and $a^{\rm d} \cdot b = b^{\rm d} \cdot
    a = 0$. In this case, all points of the lattice are projected, resulting
    in a dense set of positions on the horizontal line at the bottom. At each
    such position, a diffraction peak is located. Its height, i.e., the
    intensity seen in an experiment, is determined by the vertical coordinate
    $k_{\rm int}$ of the unique corresponding point of the dual lattice. The
    explicit value is given by the function $I(k_{\rm int})\sim \left(
      \frac{\sin (\sqrt{2}\pi k_{\rm int})}{\sqrt{2}\pi k_{\rm int}}\right)^2$
    which is displayed on the right hand side.  }
\label{diffrac}\smallskip
\end{picbox}
\end{minipage}
}}
\clearpage

\section{What are cut and project sets?}

The realization of our point sets as lingering images of lattices in higher
dimensional spaces is both visually appealing and sheds light on diffraction.
However, the use of conjugation as we used it appears as a miracle and one is
left wondering why it worked and when we might expect it to work again. In
fact, the answer to this is not really known. We do not know when a given
aperiodic point set, even if it is pure point diffractive, may be realized in
the cut and project formalism. We do know that they are not restricted to sets
involving irrationalities like $\sqrt 2$. One of the most interesting and
earliest examples of this is the one based on the Robinson square tiles.

These tiles arose out of another one of the streams whose confluence produced
the subject of aperiodic order, namely the decision problem for tilings. Given
a finite number of tile types, is there an algorithm for determining whether
or not the plane can be tiled (covered without gaps and overlaps) by
translated copies of these tiles?  This problem had been raised and later
brought to a negative conclusion by logicians.  Tiles that only can tile
aperiodically lie at the heart of this undecidability, and the hunt was on for
the smallest collections of such tiles.

Raphael Robinson made a very interesting contribution to this by first linking
the problem of tiling a plane with marked square tiles to Turing machines and
the famous Halting Problem, and also coming up with a simple set of $6$ square
tiles with markings (shown in Box~\ref{robinson} --- actually 28 tiles since
all rotated and reflected images are also to be included) that only tile
aperiodically. A rather dramatic proof of this can be glimpsed from the
subsequent pictures where it is seen that legal arrangements of the tiles lead
to a family of interlocking squares of increasing (by factors of $2$) sizes.
The aperiodicity is obvious: no finite translation could take the squares of
all sizes into themselves.

If we mark the centre of each tile by a coloured point (to indicate its type)
then we get $6$ (or $28$) families of points which are subsets of a square
lattice. These point sets are in fact cut and project sets, but now the
`higher dimensional' space is far more exotic: it is the product of a
Euclidean plane and an arithmetical-topological space that is based on the
so-called $2$-adic numbers. In spite of being very different from a Euclidean
space, the diffraction results are provable much as before. Each of these
point sets is pure point diffractive!

There remains though, the difficult problem of characterizing cut and project
sets.

\bigskip
\centerline{
\fbox{
\begin{minipage}{0.9\textwidth}
\[
\centerline{\hspace*{0.05\textwidth}%
\epsfxsize=0.232\textwidth\epsfbox{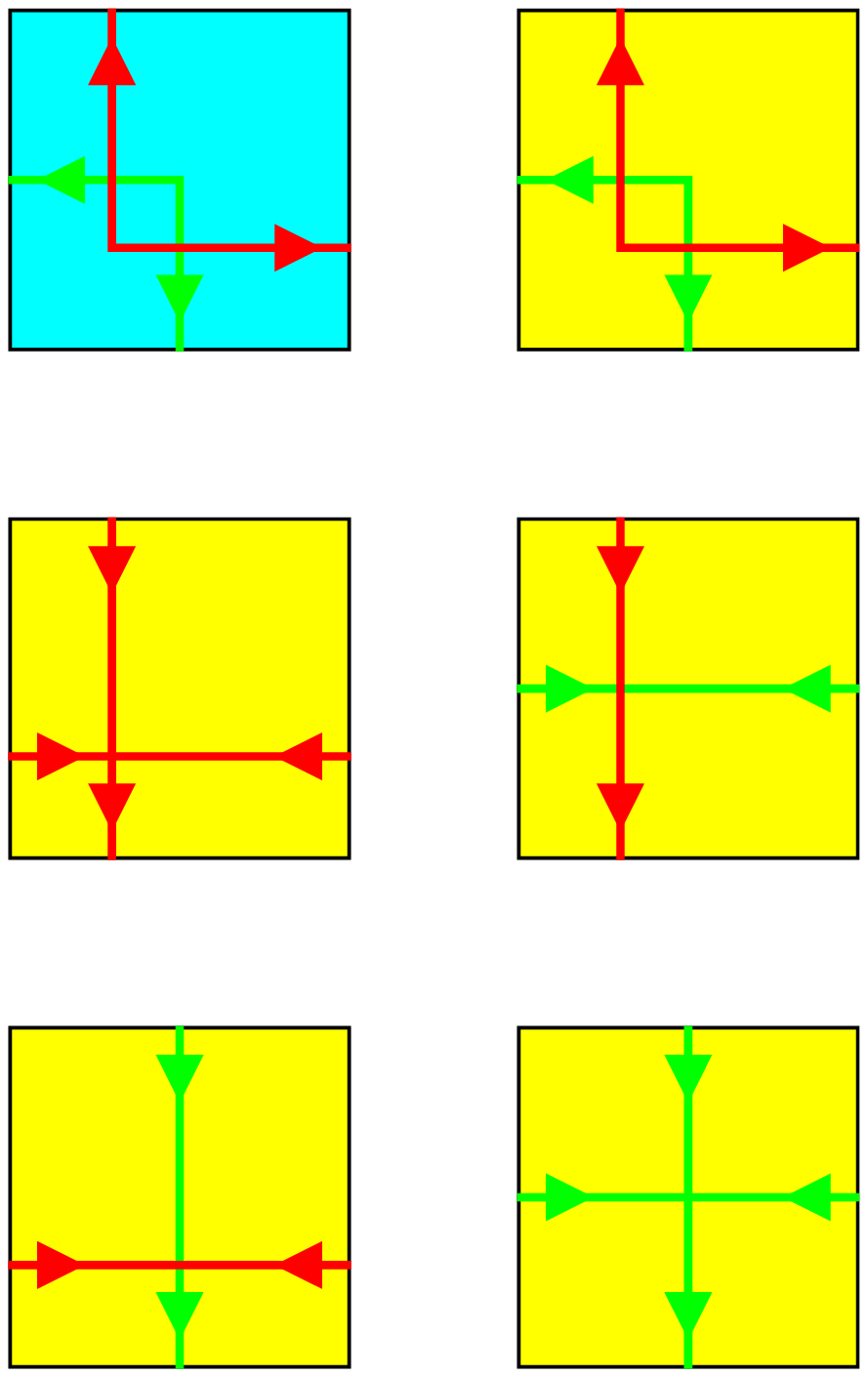}\hfill%
\epsfxsize=0.432\textwidth\epsfbox{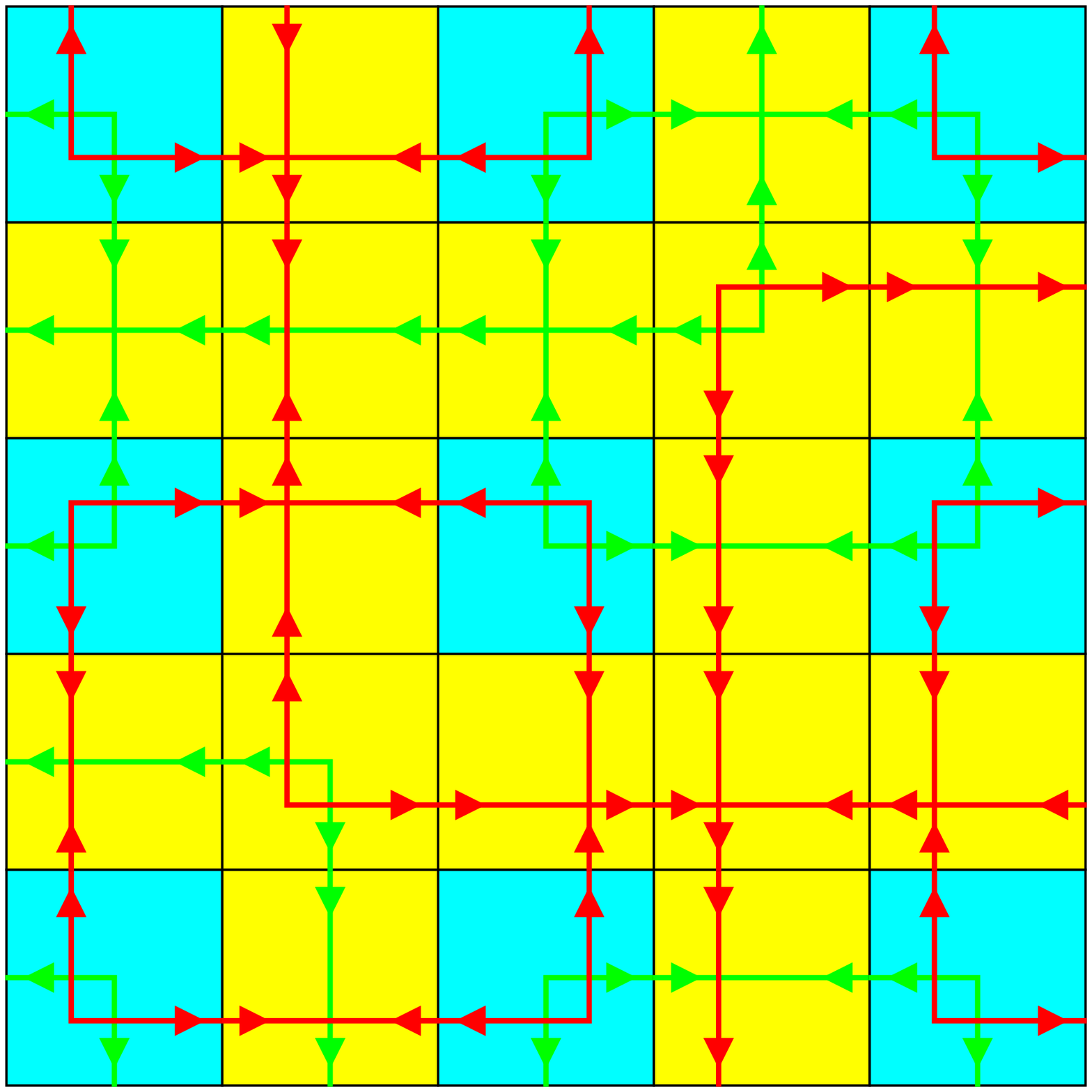}\hspace*{0.05\textwidth}}
\]
\[
\centerline{\epsfxsize=0.9\textwidth\epsfbox{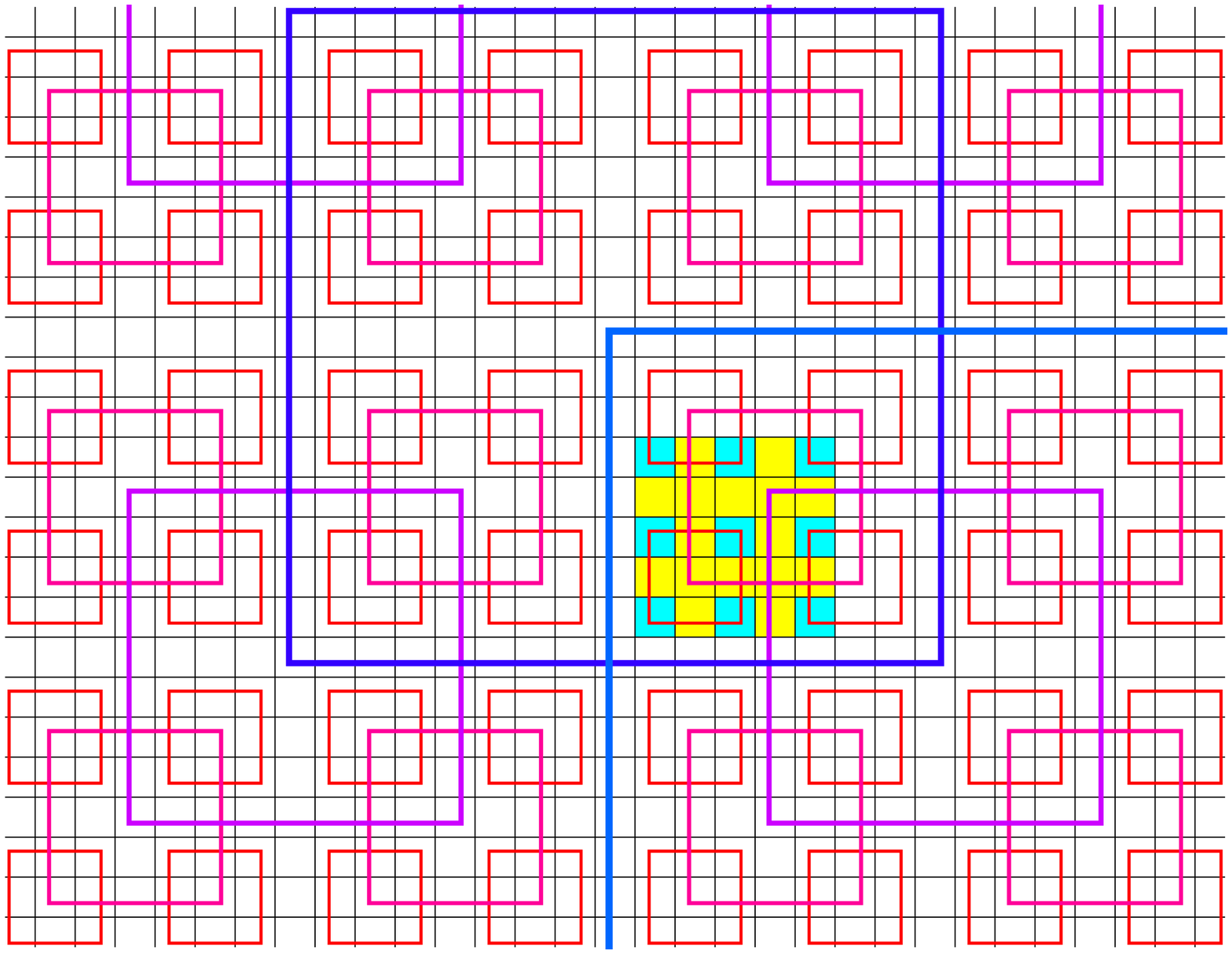}}
\]
\begin{picbox}
  Robinson tiling\newline {\em The six Robinson tiles (upper left) given as
    squares of two different colours that are labeled by two types of oriented
    lines. Together with their images under rotation and reflection they make
    up an aperiodic set of tiles, if one requires that the oriented lines
    match at the edges, and that exactly three yellow squares meet at each
    corner (upper right). Disregarding the green lines, the red lines make up
    a pattern of interlocking larger and larger squares, indicated by
    different colours in the lower picture. The region tiled by coloured
    squares corresponds to the patch shown above.  }
\label{robinson}\smallskip
\end{picbox}
\end{minipage}
}}
\bigskip

\section{Probabilistic ideas}

As was briefly mentioned in the beginning, quasicrystals can also be seen as a
stepping stone for bridging the gap between perfect crystals on the one
extreme and amorphous solids on the other. It can clearly only be a first
step, as we have seen how close they are to crystals in so many properties.

Indeed, as all constructions above have shown, quasicrystals are completely
deterministic, and what is still missing here is a source for some kind of
randomness, or stochastic disorder. This would be an entire story in itself,
but we can at least indicate one way to use crystallographic and
quasicrystallographic tilings to make some steps into this new direction.  The
new facet here is that the underlying mechanism is {\em statistical} in
origin, both for the reason of existence and for the appearance of symmetries,
which are also statistical now.

Inspecting Box \ref{abpatch} again, we now remove all markings, and also the
long edges of the triangles. We obtain a square-rhombus tiling, with many
``simpletons''. By these we mean little (irregular) hexagons built from one
square and two rhombi, as shown in Box \ref{flip}.  They can now be flipped as
indicated, without affecting any face-to-face condition. If we randomly pick
such simpletons and flip them, and continue doing so for a while (in fact, for
eternity), we arrive at what is called the square-rhombus random tiling
ensemble.  A snapshot is shown in Box \ref{randomtiling}.

In this way, we have introduced an element of randomness into our tiling, but
without destroying the basic building blocks (the square and the rhombus) and
their face-to-face arrangements. Also, this does not change the ratio of
squares to rhombi. Nevertheless, there are many such tilings now, in fact even
exponentially many, i.e.\ the number of different patches of a given size
grows exponentially with the size!  This means that the ensemble even has
positive entropy density, which opens the door for a completely different
explanation of why we see them in nature: they are, given the building blocks
(e.g.\ in the form of rather stable atomic clusters that can agglomerate),
``very likely''. Recent evidence seems to point into this direction, and a
more detailed investigation of these random tilings is desirable.

In fact, one could even start from just a pool of tiles of both types and
admit all assemblies that cover the plane without gaps or overlaps, and
without violating the face-to-face condition of the tiles.  This way, one gets
an even larger class of tilings, called the unrestricted square-rhombus random
tiling ensemble, where arbitrary ratios of squares to rhombi are realizable.
Among them, we also find the ones constructed by randomization of perfect
tilings as explained above, and one can show that the tilings of maximal
entropy (which basically means the most likely ones of this enlarged ensemble)
have the square-rhombi ratio of the perfect Ammann-Beenker pattern and show
eightfold, hence maximal, symmetry! The latter has to be interpreted in the
statistical sense, meaning that each patch one can find occurs in all 8
orientations with the same frequency.  This brings about a totally different
symmetry concept which is statistical rather than deterministic in origin, a
somewhat puzzling thought perhaps.  Nevertheless, this is sufficient to make
the corresponding diffraction image exactly eightfold symmetric!
\bigskip\bigskip\bigskip

\centerline{
\fbox{
\begin{minipage}{0.9\textwidth}
\[
\centerline{\epsfxsize=0.9\textwidth\epsfbox{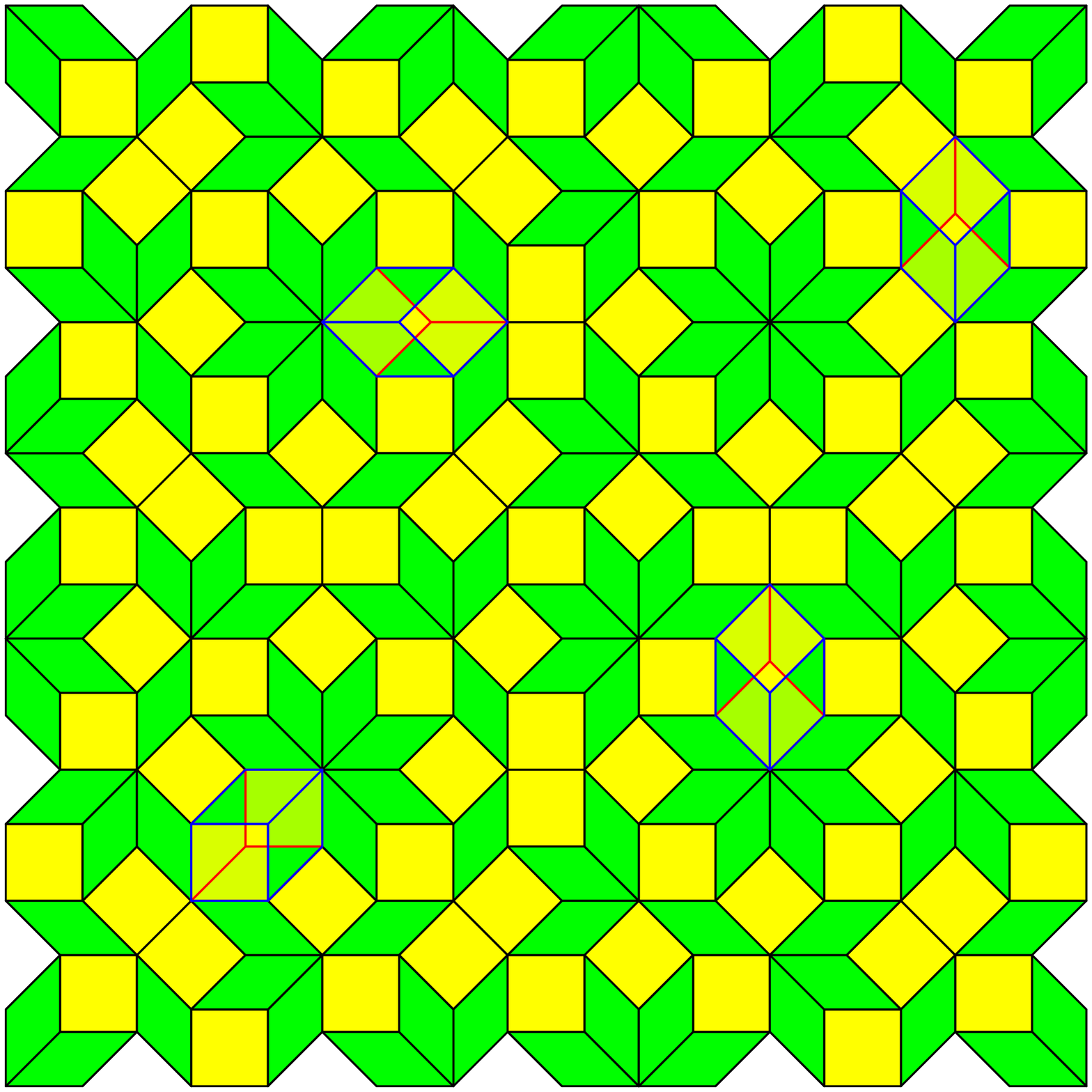}}
\]
\begin{picbox}
  Simpleton flips\newline {\em Four examples of simpleton flips in a patch of
    the perfect Ammann-Beenker tiling. The hexagons and their original
    dissection into a square and two rhombi are marked by the blue lines,
    whereas the red lines indicate the flipped arrangement. Note that only the
    three internal lines in the hexagon are affected by the flip, the outer
    shape stays the same. One can view the patch, and all variants obtained by
    such elementary simpleton flips, also as the projection of a (fairly
    rugged) roof in 3-space --- the two versions of the simpleton fillings
    then correspond to the projection of two different half surfaces of a
    cube.  }
\label{flip}\smallskip
\end{picbox}
\end{minipage}
}}
\bigskip

\bigskip
\centerline{
\fbox{
\begin{minipage}{0.9\textwidth}
\[
\centerline{\epsfxsize=\textwidth\epsfbox{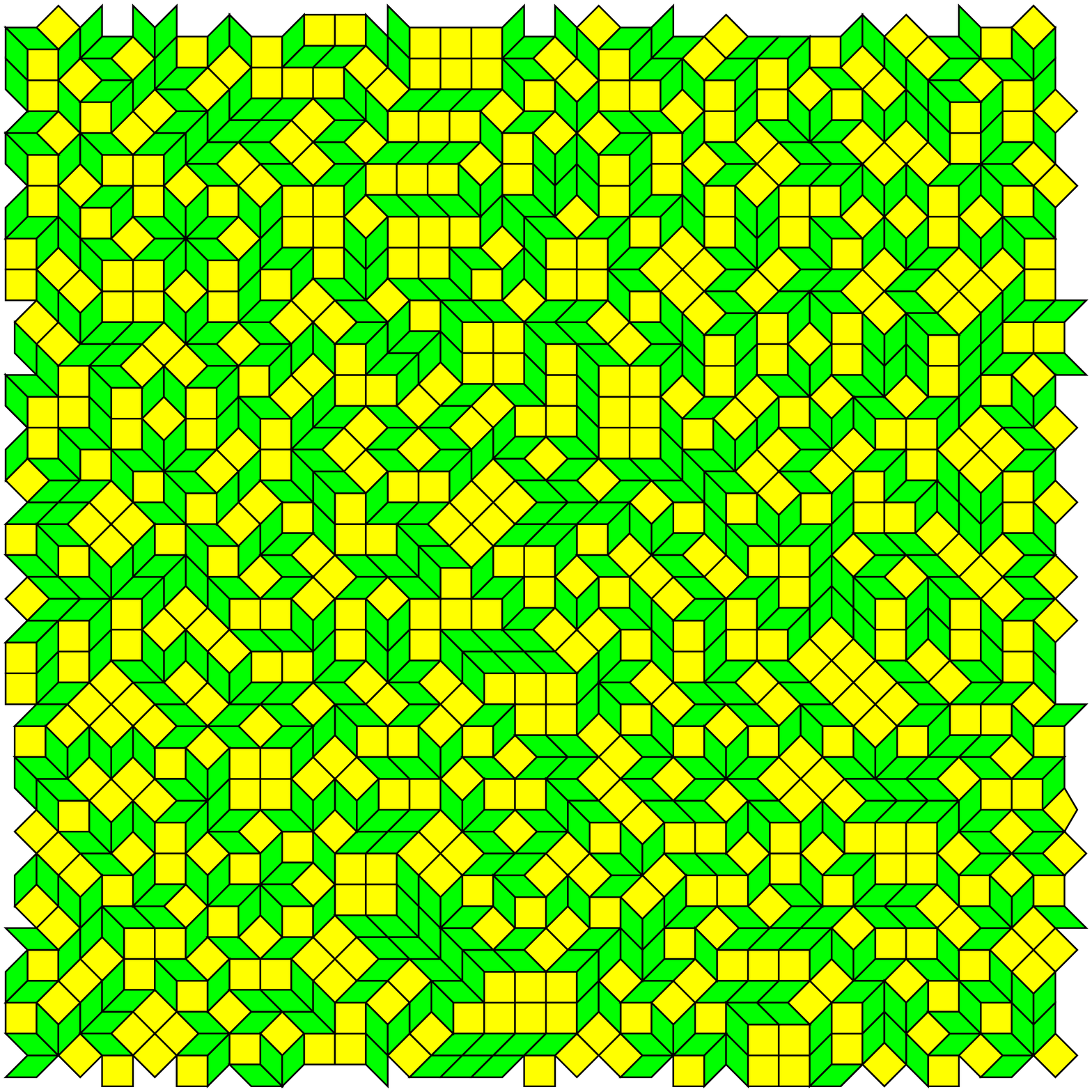}}
\]
\begin{picbox}
  Square-rhombus random tiling\newline {\em A patch of a square-rhombus random
    tiling obtained by randomly rearranging a large approximating patch of the
    perfect Ammann-Beenker tiling. In fact, we started from a square-shaped
    patch as those shown in Box \ref{sequence}, whose translated copies, when
    glued together along its boundaries, generate a periodic pattern that
    violates the perfect matching rules only in the corners where the pieces
    are glued together. The same procedure could be applied to the disordered
    patch shown here, resulting in a periodic pattern which simply has an
    enormously large building block, namely the one shown above!}
\label{randomtiling}\smallskip
\end{picbox}
\end{minipage}
}}
\bigskip

\bigskip
\centerline{
\fbox{
\begin{minipage}{0.9\textwidth}
\[
\centerline{\epsfxsize=\textwidth\epsfbox{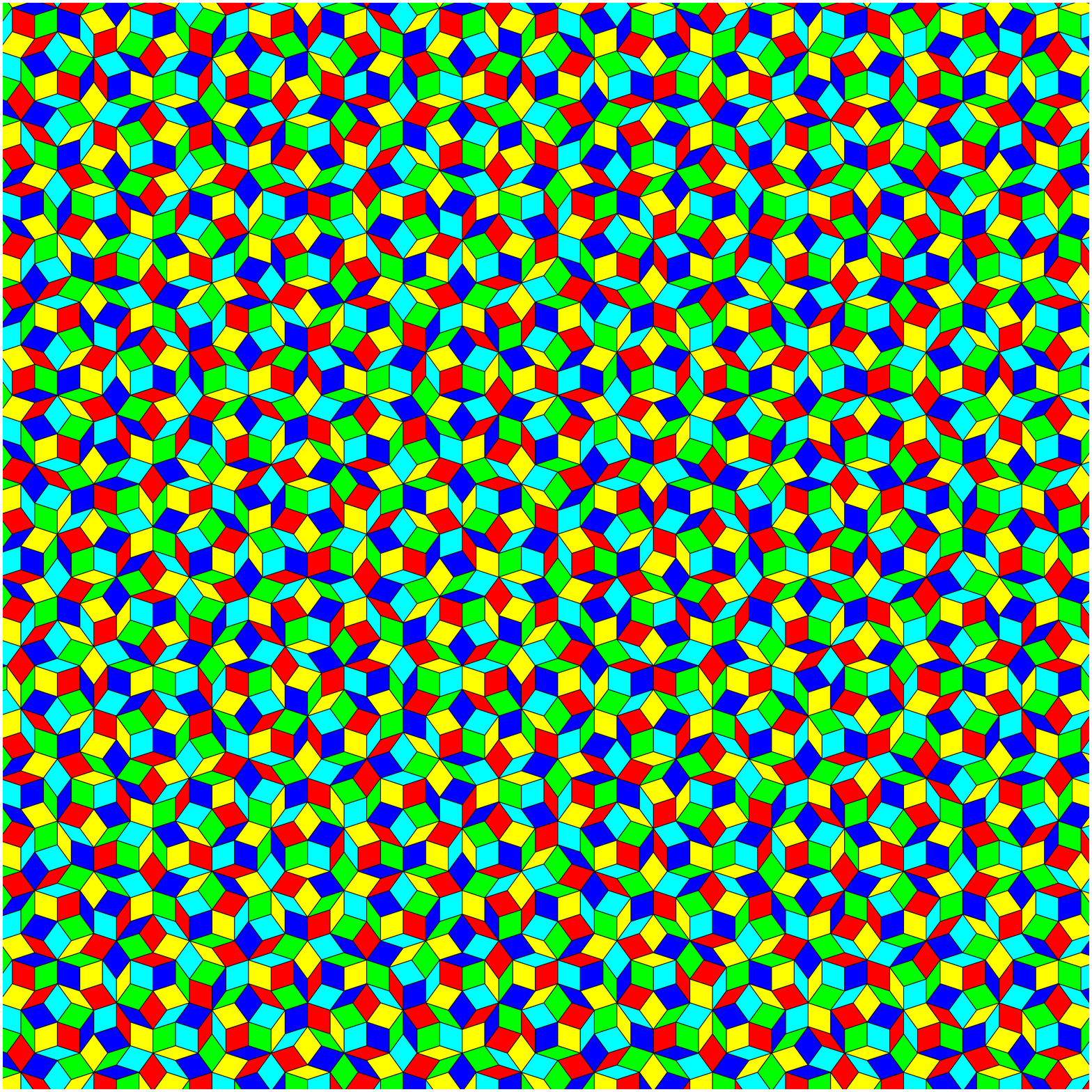}}
\]
\begin{picbox}
  A colour-symmetric Penrose tiling\newline {\em The picture shows a colouring
    of the Penrose tiling with five different colours. The colours are chosen
    such that they permute in a definite way under rotation of the tiling.
    Figure courtesy of Max Scheffer (Chemnitz).}
\label{colour}\smallskip
\end{picbox}
\end{minipage}
}}
\clearpage

\section{Summing up}

One fascinating thing about the type of order exemplified in this discussion
is how very close it comes to being periodic without admitting any actual
periods.

So, let us ask again: `what is aperiodic order?'.  At present, we have a
reasonable qualitative and a partial quantitative understanding, some aspects
of which we have tried to explain above. However, we still don't have a
complete answer, and such an answer might lie well into the future.

But what we do know is that there is a universe of beautiful questions out
there, with unexpected results to be found, and with many cross-connections
between seemingly disjoint disciplines.  On top of that, it is definitely a
lot of fun, for example, when producing new variants of Penrose tilings with
colour symmetries, such as the example shown in Box \ref{colour} below!  For a
recent bibliographical review of the literature, we refer the reader to
\cite{BM}.

\end{document}